\documentclass[journal]{new-aiaa}
\usepackage[utf8]{inputenc}
\usepackage{textcomp}

\usepackage{graphicx}

\usepackage{amsmath,amssymb,amsfonts}
\allowdisplaybreaks
\usepackage{bm}
\usepackage{relsize}
\usepackage{hhline}
\usepackage{ulem}
\usepackage[version=4]{mhchem}
\usepackage{siunitx}
\usepackage{longtable,tabularx}
\usepackage{color}
\usepackage{float}
\usepackage{adjustbox}
\usepackage{subfig}
\usepackage{mathtools}
\usepackage{units}
\usepackage{epstopdf}

\definecolor{commentgreen}{RGB}{0, 153, 0}
\definecolor{highlightred}{RGB}{153,0,0}

\newcommand{\sref}[1]{\S\ref{#1}}
\newcommand{\LOS}{{LoS }}

\newcommand{\real}{{\mathbb{R}}}
\newcommand{\realp}{\real_+}
\newcommand{\realpp}{\real_{++}}

\newcommand{\ind}{i}

\newcommand{\der}[2]{\frac{d #1}{d #2}}

\renewcommand{\unit}[1]{\mathrm{#1}}

\newcommand{\inertial}{\mathcal{I}}
\newcommand{\body}{\mathcal{B}}

\newcommand{\cframe}[1]{\mathcal{F}_{#1}}
\newcommand{\Fbody}{\cframe{\body}}
\newcommand{\Finertial}{\cframe{\inertial}}

\newcommand{\xI}{\bm{x}_{\inertial}}
\newcommand{\yI}{\bm{y}_{\inertial}}
\newcommand{\zI}{\bm{z}_{\inertial}}
\newcommand{\xB}{\bm{x}_{\body}}
\newcommand{\yB}{\bm{y}_{\body}}
\newcommand{\zB}{\bm{z}_{\body}}

\newcommand{\eps}{\epsilon}
\newcommand{\dg}{^{\circ}}

\newcommand{\diag}[1]{\mathop{\bf diag}\left\lbrace #1 \right\rbrace}
\newcommand{\eye}[1]{I_{#1}}
\newcommand{\zeros}[2]{0_{#1 \times #2}}
\newcommand{\ones}[1]{\bm{1}_{#1}}
\newcommand{\dom}[1]{\mathop{\bf dom} #1}

\newcommand{\definedas}{\coloneqq}
\newcommand{\tran}{^{\top}}
\newcommand{\ipQ}{\boldsymbol{\cdot}}
\newcommand{\ipDQ}{\mathrel{\raisebox{1pt}{\scalebox{0.6}{$\circ$}}}}

\newcommand{\tf}{t_f}
\newcommand{\ti}{t_0}
\newcommand{\dtau}{\Delta \tau}

\newcommand{\gI}{\bm{g}_{\inertial}}
\newcommand{\gB}{\bm{g}_{\body}}

\newcommand{\dqgB}{\dq{g}_{\body}}
\newcommand{\J}{\bm{J}}

\newcommand{\m}{m}

\newcommand{\mi}{\m_0}

\newcommand{\rI}{\bm{r}_{\inertial}}

\newcommand{\rB}{\bm{r}_{\body}}

\newcommand{\rIdot}{\dot{\bm{r}}_{\inertial}}

\newcommand{\ri}[1]{\bm{r}_{#1,0}}
\newcommand{\rf}[1]{\bm{r}_{#1,f}}

\newcommand{\vI}{\bm{v}_{\inertial}}
\newcommand{\vB}{\bm{v}_{\body}}
\newcommand{\vBdot}{\dot{\bm{v}}_{\body}}

\newcommand{\vi}[1]{\bm{v}_{#1,0}}
\newcommand{\vf}[1]{\bm{v}_{#1,f}}

\newcommand{\wB}{\bm{\omega}_{\body}}
\newcommand{\wBdot}{\dot{\bm{\omega}}_{\body}}

\newcommand{\wi}{\bm{\omega}_{\body,0}}
\newcommand{\wf}{\bm{\omega}_{\body,f}}

\newcommand{\FB}{\bm{F}_{\body}}
\newcommand{\TB}{\bm{M}_{\body}}
\newcommand{\uB}{\bm{u}_{\body}}

\newcommand{\duB}{\dot{\bm{u}}_{\body}}

\newcommand{\uk}[1]{\bm{u}_{#1,\ind}}

\newcommand{\ru}{\bm{r}_{u}}

\newcommand{\mwet}{m_{\text{wet}}}
\newcommand{\mdry}{m_{\text{dry}}}

\newcommand{\wmax}{\omega_{\max}}
\newcommand{\gimbal}{\delta}
\newcommand{\gimbalmax}{\delta_{\max}}
\newcommand{\los}{\xi}
\newcommand{\losmax}{\los_{\max}}

\newcommand{\losvec}{\bm{p}_{\body}}

\newcommand{\gs}{\gamma}
\newcommand{\gsmax}{\gamma_{\max}}

\newcommand{\gsf}{c_g}
\newcommand{\gsM}{M_g}
\newcommand{\tiltmax}{\theta_{\max}}
\newcommand{\tilt}{\theta}
\newcommand{\tiltf}{c_t}
\newcommand{\tiltM}{M_t}

\newcommand{\umin}{u_{\min}}
\newcommand{\umax}{u_{\max}}

\newcommand{\dgimmax}{\dot{\gimbal}_{\max}}

\newcommand{\stcz}{\bm{z}}
\newcommand{\stcnz}{n_z}
\newcommand{\stcng}{n_g}

\newcommand{\stcg}{g}
\newcommand{\stcc}{c}
\newcommand{\stch}{h}
\newcommand{\shat}{\sigma}

\newcommand{\setq}{\mathbb{Q}}
\newcommand{\setuq}{\setq_u}

\newcommand{\setdq}{\widetilde{\setq}}
\newcommand{\setudq}{\widetilde{\setq}_u}
\newcommand{\setuqR}{\real_{u}^4}
\newcommand{\setudqR}{\real_{u}^8}
\newcommand{\quat}[1]{\bm{#1}}
\newcommand{\q}{\bm{q}}
\newcommand{\qv}{\q_v}
\newcommand{\qs}{q_4}

\newcommand{\qid}{\q_{\text{id}}}
\newcommand{\qf}{\q_f}
\newcommand{\qb}{\bar{\q}}
\newcommand{\qdot}{\dot{\q}}
\newcommand{\dq}[1]{\tilde{\bm{#1}}}
\newcommand{\dqq}{\dq{q}}
\newcommand{\dqqdot}{\dot{\dqq}}
\newcommand{\dqb}{\bar{\dqq}}
\newcommand{\dqw}{\dq{\omega}}
\newcommand{\dqwdot}{\dot{\dqw}}

\newcommand{\qskew}[1]{[ #1 ]_{\otimes}}
\newcommand{\qskewstar}[1]{[ #1 ]_{\otimes}^*}

\newcommand{\dqqig}{\bm{b}_{\dqq}}
\newcommand{\dqwig}{\bm{b}_{\dqw}}

\newcommand{\Jtr}{J_{\text{tr}}}
\newcommand{\Jvc}{J_{\text{vc}}}
\newcommand{\xxb}{\bar{\bm{x}}}
\newcommand{\uub}{\bar{\bm{u}}}
\newcommand{\ssb}{\bar{s}}
\newcommand{\zzb}{\bar{\bm{z}}}
\newcommand{\zz}{\bm{z}}
\newcommand{\xx}{\bm{x}}
\newcommand{\XX}{\bm{X}}
\newcommand{\uu}{\bm{u}}
\newcommand{\UU}{\bm{U}}
\newcommand{\ww}{R}
\renewcommand{\ss}{s}
\renewcommand{\SS}{S}
\newcommand{\dummy}{\zeta}
\newcommand{\ee}{\bm{\eta}}
\newcommand{\eek}{\eta_{\ind}}
\newcommand{\vv}{\bm{\nu}}
\newcommand{\VV}{\bm{V}}
\newcommand{\wtr}{\bm{w}_{\text{tr}}}
\newcommand{\wtrk}{\bm{w}_{\text{tr},\ind}}
\newcommand{\wvc}{w_{\text{vc}}}
\newcommand{\res}{\Delta}
\newcommand{\resk}{\res_\ind}

\newcommand{\resmin}{\res_{\min}}
\newcommand{\setK}{\mathsf{N}}
\newcommand{\setKb}{\bar{\mathsf{N}}}
\newcommand{\dxtol}{\delta x_{\text{tol}}}

\newcommand{\Px}{P_x}
\newcommand{\Pu}{P_u}
\newcommand{\qx}{\bm{p}_x}
\newcommand{\qu}{\bm{p}_u}
\newcommand{\Pt}{p_t}

\newcommand{\xxs}{\hat{\xx}}

\newcommand{\uus}{\hat{\uu}}
\newcommand{\sss}{\hat{\ss}}
\newcommand{\ees}{\hat{\ee}}
\newcommand{\eeks}{\hat{\eta}_{\ind}}
\newcommand{\XXs}{\hat{\XX}}
\newcommand{\UUs}{\hat{\UU}}
\newcommand{\VVs}{\hat{\VV}}

\newtheorem{theorem}{Theorem}[section]

\newtheorem{lemma}[theorem]{Lemma}

\newtheorem{remark}[theorem]{Remark}

\newtheorem{problem}{Problem}

\newcommand*\samethanks[1][\value{footnote}]{\footnotemark[#1]}

\usepackage{tikz}
\usetikzlibrary{math,calc,patterns,decorations.pathmorphing,decorations.markings}
\usepackage{pgfplots}
\pgfplotsset{compat=1.12}
\usepackage{ifthen}
\usepackage{rotating}
\usetikzlibrary{patterns}
\usetikzlibrary{math}
\usetikzlibrary{scopes}
\usetikzlibrary{fadings}
\usetikzlibrary{arrows}
\tikzset{>=latex}

\pgfdeclarelayer{bg}    
\pgfsetlayers{bg,main}  

\definecolor{beige}{RGB}{245,245,220}
\definecolor{darkred}{rgb}{0.90,0.00,0.00}
\definecolor{darkgreen}{rgb}{0.00,0.45,0.00}
\definecolor{midgreen}{rgb}{0.00,0.65,0.00}
\definecolor{purple}{rgb}{0.50,0.00,1.00}

\definecolor{ucol}{RGB}{255,0,0}
\definecolor{gcol}{RGB}{0,120,0}
\definecolor{scol}{RGB}{63, 226, 45}
\definecolor{vcol}{RGB}{0,0,0}
\definecolor{acol}{RGB}{0, 128, 255}
\definecolor{dcol}{RGB}{204,102,0}
\definecolor{lcol}{RGB}{204,102,0}
\definecolor{bcol}{RGB}{0,0,0}
\definecolor{ocol}{RGB}{167,167,167}

\definecolor{propcolor}{rgb}{1.0,0.0,0.0}
\definecolor{solvecolor}{rgb}{0.0,0.55,0.0}
\definecolor{scalecolor}{RGB}{75, 47, 132}
\definecolor{fillcolor123}{rgb}{1,1,1}

\tikzstyle{block} = [draw, ultra thick, fill=blue!20, rectangle, 
    rounded corners=6pt, minimum height=3em, minimum width=6em]
\tikzstyle{mexblock} = [draw, ultra thick, magenta, fill=blue!20, rectangle, 
    rounded corners=6pt, minimum height=3em, minimum width=6em, text=black]
\tikzstyle{sum} = [draw, fill=blue!20, circle, node distance=2cm]
\tikzstyle{input} = [coordinate]
\tikzstyle{output} = [coordinate]
\tikzstyle{pinstyle} = [pin edge={to-,thin,black}]

\newcommand{\lwarrow}{0.50mm}
\newcommand{\lwthin}{0.10mm}
\newcommand{\lwthick}{0.50mm}
\newcommand{\lwtext}{0.25mm}

\newcommand{\threeaxes}[8]{
	\tikzmath{
    	\Lxr=#3;\Lxl=#3;\Lyt=#4;\Lyb=#4;\Zt= #5;\Zb= #6;\Xang=#7;\Yang=#8;
    	\Xxr= cos(\Xang)*\Lxr; \Xyr=-sin(\Xang)*\Lxr;
    	\Xxl=-cos(\Xang)*\Lxl; \Xyl= sin(\Xang)*\Lxl;
        \Yxt= sin(\Yang)*\Lyt; \Yyt= cos(\Yang)*\Lyt;
        \Yxb=-sin(\Yang)*\Lyb; \Yyb=-cos(\Yang)*\Lyb;
        \zzz=0;
    }
    \begin{scope}[shift={(#1,#2)},rotate=0]
        \ifx\Zt\zzz\else  \draw[black,->] (0,0) -- +(0, \Zt);    \fi
        \ifx\Zb\zzz\else  \draw[black,->] (0,0) -- +(0,-\Zb);    \fi
        \ifx\Lxr\zzz\else \draw[black,->] (0,0) -- +(\Xxr,\Xyr); \fi
        \ifx\Lxl\zzz\else \draw[black,->] (0,0) -- +(\Xxl,\Xyl); \fi
        \ifx\Lyt\zzz\else \draw[black,->] (0,0) -- +(\Yxt,\Yyt); \fi
        \ifx\Lyb\zzz\else \draw[black,->] (0,0) -- +(\Yxb,\Yyb); \fi
	\end{scope}
}

\newcommand{\isoaxes}[9]{
    \tikzmath{
        \rot=#3;
        \len=#4;
        \ddd=#5;
        \pX=  0; \Xx=cos(\pX)*\len; \Xy=sin(\pX))*\len;
        \pY=120; \Yx=cos(\pY)*\len; \Yy=sin(\pY))*\len;
        \pZ=240; \Zx=cos(\pZ)*\len; \Zy=sin(\pZ))*\len;
        \Xxx=cos(\ddd)*\Xx-sin(\ddd)*\Xy; \Xxy=sin(\ddd)*\Xx+cos(\ddd)*\Xy;
        \Yxx=cos(\ddd)*\Yx-sin(\ddd)*\Yy; \Yxy=sin(\ddd)*\Yx+cos(\ddd)*\Yy;
        \Zxx=cos(\ddd)*\Zx-sin(\ddd)*\Zy; \Zxy=sin(\ddd)*\Zx+cos(\ddd)*\Zy;
    }
    \begin{scope}[shift={(#1,#2)},rotate=\rot]
		\filldraw[black] (0,0) circle (2pt);
        \draw[black,thick,->] (0,0) -- +(\Xx,\Xy);
        \draw[black,thick,->] (0,0) -- +(\Yx,\Yy);
        \draw[black,thick,->] (0,0) -- +(\Zx,\Zy);
        \draw (0.1,0.6)   node[rotate=0,anchor=center] {#6};
        \draw (\Xxx,\Xxy) node[rotate=0,anchor=center] {#7};
        \draw (\Yxx,\Yxy) node[rotate=0,anchor=center] {#8};
        \draw (\Zxx,\Zxy) node[rotate=0,anchor=center] {#9};
    \end{scope}
}

\newcommand{\isoaxesNL}[6]{
    \tikzmath{
        \rot=#3;
        \len=#4;
        \ddd=#5;
        \pX=  0; \Xx=cos(\pX)*\len; \Xy=sin(\pX))*\len;
        \pY=120; \Yx=cos(\pY)*\len; \Yy=sin(\pY))*\len;
        \pZ=240; \Zx=cos(\pZ)*\len; \Zy=sin(\pZ))*\len;
        \Xxx=cos(\ddd)*\Xx-sin(\ddd)*\Xy; \Xxy=sin(\ddd)*\Xx+cos(\ddd)*\Xy;
        \Yxx=cos(\ddd)*\Yx-sin(\ddd)*\Yy; \Yxy=sin(\ddd)*\Yx+cos(\ddd)*\Yy;
        \Zxx=cos(\ddd)*\Zx-sin(\ddd)*\Zy; \Zxy=sin(\ddd)*\Zx+cos(\ddd)*\Zy;
    }
    \begin{scope}[shift={(#1,#2)},rotate=\rot]
		\filldraw[black] (0,0) circle (2pt);
        \draw[black,thick,->] (0,0) -- +(\Xx,\Xy);
        \draw[black,thick,->] (0,0) -- +(\Yx,\Yy);
        \draw[black,thick,->] (0,0) -- +(\Zx,\Zy);
        \draw (0.1,0.6)   node[rotate=0,anchor=center] {#6};
    \end{scope}
}

\newcommand{\isoaxesNLC}[7]{
    \tikzmath{
        \rot=#3;
        \len=#4;
        \ddd=#5;
        \pX=  0; \Xx=cos(\pX)*\len; \Xy=sin(\pX))*\len;
        \pY=120; \Yx=cos(\pY)*\len; \Yy=sin(\pY))*\len;
        \pZ=240; \Zx=cos(\pZ)*\len; \Zy=sin(\pZ))*\len;
        \Xxx=cos(\ddd)*\Xx-sin(\ddd)*\Xy; \Xxy=sin(\ddd)*\Xx+cos(\ddd)*\Xy;
        \Yxx=cos(\ddd)*\Yx-sin(\ddd)*\Yy; \Yxy=sin(\ddd)*\Yx+cos(\ddd)*\Yy;
        \Zxx=cos(\ddd)*\Zx-sin(\ddd)*\Zy; \Zxy=sin(\ddd)*\Zx+cos(\ddd)*\Zy;
    }
    \begin{scope}[shift={(#1,#2)},rotate=\rot]
		\filldraw[#7] (0,0) circle (2pt);
        \draw[#7,thick,->] (0,0) -- +(\Xx,\Xy);
        \draw[#7,thick,->] (0,0) -- +(\Yx,\Yy);
        \draw[#7,thick,->] (0,0) -- +(\Zx,\Zy);
        \draw (0.1,0.6)   node[#7,rotate=0,anchor=center] {#6};
    \end{scope}
}

\newcommand{\pane}[7]{
	\tikzmath{
    	\rot=#3;
    	\width=#4;
        \height=#5;
        \corner=#6;
    	\px1= 0.5*\width-\corner; \py1=-0.5*\height;
        \px2= 0.5*\width;         \py2=-0.5*\height+\corner;
        \px3= 0.5*\width;         \py3= 0.5*\height-\corner;
        \px4= 0.5*\width-\corner; \py4= 0.5*\height;
        \px5=-0.5*\width+\corner; \py5= 0.5*\height;
        \px6=-0.5*\width;         \py6= 0.5*\height-\corner;
        \px7=-0.5*\width;         \py7=-0.5*\height+\corner;
        \px8=-0.5*\width+\corner; \py8=-0.5*\height;
    }
	\begin{scope}[shift={(#1,#2)},rotate=\rot]
    	\filldraw[{#7}]
        (\px1,\py1) to[out=    0,in=  -90] (\px2,\py2) --
        (\px3,\py3) to[out=   90,in=    0] (\px4,\py4) --
        (\px5,\py5) to[out= -180,in=   90] (\px6,\py6) --
        (\px7,\py7) to[out=  -90,in= -180] (\px8,\py8) -- cycle;
    \end{scope}
}
\newcommand{\cpane}[7]{
	\tikzmath{
		\pLTx=#1;
		\pRBx=#2;
		\pLTy=#3;
		\pRBy=#4;
	}
	\pane{0.5*\pLTx+0.5*\pRBx}{0.5*\pLTy+0.5*\pRBy}{#5}{\pRBx-\pLTx}{\pLTy-\pRBy}{#6}{#7}
}

\newcommand{\coneback}[7]{
	\tikzmath{\rot=#3;
              \length=#4; \radius=\length*tan(0.5*#5); \depth=#6;
              \sx =  cos(\rot)*#1 + sin(\rot)*#2;
              \sy = -sin(\rot)*#1 + cos(\rot)*#2;
    }
    \begin{scope}[shift={(\sx,\sy)},transform canvas={rotate=\rot}]
	    \draw[{#7}] (\radius,-\length) arc(360:180: {\radius} and {-\depth});
    \end{scope}
}

\newcommand{\cone}[7]{
	\tikzmath{\rot=#3;
              \length=#4; \radius=\length*tan(0.5*#5); \depth=#6;
              \sx =  cos(\rot)*#1 + sin(\rot)*#2;
              \sy = -sin(\rot)*#1 + cos(\rot)*#2;
	}
    \begin{scope}[shift={(\sx,\sy)},transform canvas={rotate=\rot}]
    	\fill[{#7}] (0,0) -- (\radius,-\length) arc(360:180: {\radius} and {\depth}) -- cycle;
    	\draw[color=black!100] (0,0) -- (\radius,-\length) arc(360:180: {\radius} and {\depth}) -- cycle;
    \end{scope}
}

\newcommand{\rocket}[6]{
	\tikzmath{\rot=#3;
    		  \length=#4;
              \throttle=(\length/0.8)*0.6*#5;
              \gimbalangle=#6;
    		  \LL = \length; \RR = \LL/8;
              \HH = \LL/5;   \WW = \LL/16;
              \rr = \LL/16;
              \HG = \LL/8;   \RG = \LL/11; \WG = \LL/16;
              \sx =  cos(\rot)*#1 + sin(\rot)*#2;
              \sy = -sin(\rot)*#1 + cos(\rot)*#2;
    }
    
    \begin{scope}[shift={(\sx,\sy)},transform canvas={rotate=\rot}]
		\begin{scope}[shift={(0,-0.5*\LL)},rotate=\gimbalangle]
          	\fill[fill=orange!100,shading=axis,shading angle=90,left color=orange!100,right color=orange!25]
            	(-\RG,-\HG) -- (0,-\HG-\throttle) -- (\RG,-\HG) -- cycle;
        	\filldraw[color=black!100,fill=black!10,shading=axis,shading angle=90,left color=black!30,right color=black!0]
        	(0,0) -- (\RG,-\HG) arc(360:180: {\RG} and {\WG}) -- cycle;
		\end{scope}

    	\filldraw[color=black!100,fill=black!10,shading=ball,shading angle=90,left color=black!30,right color=black!0]
        	(-\RR,-0.5*\LL) arc(180:360: {\RR} and {\WW}) -- (\RR,0.5*\LL) -- (\RR,0.5*\LL) arc(360:180: {\RR} and {\WW}) -- (-\RR,-0.5*\LL) -- cycle;
        
    	\filldraw[color=black!100,fill=black!10,shading=axis,shading angle=90,left color=black!30,right color=black!0]
        	(\RR,0.5*\LL) arc(360:180: {\RR} and {\WW}) -- (-\rr,0.5*\LL+\HH) arc(-180:0: {\rr} and {-\WW}) -- cycle;
    \end{scope}
}

\title{Dual Quaternion Based Powered Descent Guidance \\ with State-Triggered Constraints}

\author{Taylor P. Reynolds \footnote{Ph.D. Candidate, UW Aero. and Astro., AIAA Student Member, \texttt{\{tpr,mszmuk,danylo\}@uw.edu}.}, Michael Szmuk\samethanks[1], Danylo Malyuta\samethanks[1], Mehran Mesbahi\footnote{Professor, UW Aero. and Astro., AIAA Associate Fellow, \texttt{\{mesbahi,behcet\}@uw.edu}.} and Beh\c{c}et A\c {c}{\i}kme\c{s}e\samethanks[2]}
\affil{University of Washington, Seattle, WA, 98195}
\author{John M. Carson III\footnote{SPLICE Project Manager, Johnson Space Center, AIAA Associate Fellow, \texttt{john.m.carson@nasa.gov}}}
\affil{NASA Johnson Space Center, Houston, TX, 77058}

\begin{document}

\maketitle

\begin{abstract}
    This paper presents a numerical algorithm for computing 6-degree-of-freedom free-final-time powered descent guidance trajectories. The trajectory generation problem is formulated using a unit dual quaternion representation of the rigid body dynamics, and several standard path constraints. Our formulation also includes a special line of sight constraints that is enforced only within a specified band of slant ranges relative to the landing site, a novel feature that is especially relevant to Terrain and Hazard Relative Navigation. We use the newly introduced state-triggered constraints to formulate these range constraints in a manner that is amenable to real-time implementations.  The resulting non-convex optimal control problem is solved iteratively as a sequence of convex second-order cone programs that locally approximate the non-convex problem. Each second-order cone program is solved using a customizable interior point method solver. Also introduced are a scaling method and a new heuristic technique that guide the convergence process towards dynamic feasibility. To demonstrate the capabilities of our algorithm, two numerical case studies are presented. The first studies the effect of including a slant-range-triggered line of sight constraint on the resulting trajectories. The second study performs a Monte Carlo analysis to assess the algorithm's robustness to initial conditions and real-time performance.
\end{abstract}

\section*{Nomenclature}

{\renewcommand\arraystretch{1.0}
\noindent\begin{longtable*}{@{}l @{\quad=\quad} l@{}}
$\cframe{}$ & a coordinate frame \\
$g_e$ & $9.806~\unit{m/s^2}$, standard Earth gravitational acceleration \\
$\gI$ & inertial lunar gravitational acceleration vector \\ 
$I_{\text{sp}}$ & vacuum specific impulse \\
$J$ & inertia matrix in the body coordinate frame \\
$m$ & mass \\
$\mdry$ & dry mass (i.e., zero fuel mass) \\
$N$ & number of discrete time points  \\ 
$\setK,\,\setKb$ & the sets $\{1,\ldots,N\}$ and $\{1,\ldots,N-1\}$ \\ 
$\losvec$ & boresight vector of an optical sensor in the body coordinate frame \\ 
$\Px,\qx,\Pu,\qu,\Pt$ & scaling terms associated with the state, control and time, respectively \\  
$\setq$ & quaternion manifold \\
$\setdq$ & dual quaternion manifold \\ 
$\qid$ & identity quaternion $(\zeros{3}{1},1)$ \\
$\q$ & unit quaternion representing the attitude of a rigid body \\
$\dqq$ & unit dual quaternion representing the attitude and position of a rigid body \\
$\bm{r}$ & position vector of a rigid body \\
$\real,\,\realp,\,\realpp$ & sets of real, nonnegative real, and positive real numbers \\
$\ss$ & time dilation factor relating time to scaled time \\
$t$ & time in seconds \\ 
$\umin,\,\umax$ & minimum and maximum thrust magnitude \\
$\dot{u}_{z,\max}$ & maximum throttle rate \\
$\bm{v}$ & linear velocity vector of a rigid body  \\ 
$\gs,\,\gsmax$ & current and maximum approach cone angles \\
$\gimbal,\,\gimbalmax$ & current and maximum gimbal angles \\
$\dgimmax$ & maximum gimbal rate \\
$\dxtol$ & state-based threshold used to define convergence \\ 
$\ee$ & vector of trust region radii \\ 
$\tilt,\,\tiltmax$ & current and maximum tilt angles \\
$\vv$ & virtual control \\ 
$\los,\,\losmax$ & current and maximum line of sight angles \\ 
$\sum \FB$ & sum of all externally applied forces \\
$\sum \TB$ & sum of all externally applied torques \\ 
$\shat$ & slack variable associated with state-triggered constraints \\
$\tau$ & scaled time from $0$ to $1$ \\
$\bm{\omega}$ & angular velocity of a rigid body in the body coordinate frame \\
$\wmax$ & maximum angular velocity about any axis \\
\multicolumn{2}{@{}l}{Subscripts}\\
$\body$ & resolved in the body coordinate frame \\
$f$ & final time \\ 
$\ind$ & discrete time index with respect to scaled time \\
$\inertial$ & resolved in the inertial coordinate frame \\
$0$ & initial time \\ 
\multicolumn{2}{@{}l}{Other}\\
$\bar{\cdot}$ & quantity from the previous iterate of the algorithm \\
$\hat{\cdot}$ & quantity that has been scaled \\
$\eye{n}$ & $n\times n$ identity matrix \\
$\zeros{n}{m}$ & $n\times m$ matrix of zeros \\
$\| \cdot \|_2$ & Euclidean two-norm of a vector \\
$\cdot^{\times}$ & Skew-symmetric matrix representing the cross-product operator \\
$\ipDQ$ & Dual quaternion dot product
\end{longtable*}}


\section{Introduction}\label{sec:introduction}

\lettrine{P}owered descent guidance refers to the problem of transferring a vehicle from an estimated initial state to a target state using rocket-powered engines and/or reaction control systems. While each celestial body targeted for landing presents unique design specifications, powered descent guidance methods are universally subject to constraints relating to safety, navigation, or physical vehicle requirements. Some of these constraints inherently couple the translational and rotational motion of the vehicle. This paper presents an algorithm that is designed to compute fuel-optimal trajectories in real-time that explicitly account for such constraints. Our approach is to view powered descent guidance as a trajectory generation problem in both position and attitude space for a rigid body, resulting in a 6-degree-of-freedom (6-DoF) problem. Since real-time space applications require solution times on the order of one second or less, we adopt a convex optimization based method to compute feasible and (locally) optimal trajectories~\cite{Scharf2017}. The constraint satisfaction provided by this approach permits the selection of less accessible, but more scientifically interesting, landing sites and enhances the lander’s ability to handle uncertainties~\cite{Blackmore2016,Carson2019}. Specifically,~\cite{Carson2019} provides context for this work in future descent and landing systems.

The real-time 6-DoF constrained powered descent guidance problem is difficult for several reasons. First, the nonlinear equations of motion and the presence of numerous constraints prevent both an analytical understanding of optimal solutions and the use of traditional convex optimization solution methods. Second, robust numerical methods must be developed that are insensitive to changes in problem parameters and initial solution guesses. This places an emphasis on numerical scaling, the solver that is used, and the accuracy of any approximations of the nonlinear dynamics and non-convex constraints. Third, nonlinear programming methods that achieve the two previous objectives often require computation time that exceeds what can be considered ``real-time'' for a space application. The present study addresses these challenges and provides rationale for the combination of design choices that lead to the final solution strategy.

\subsection{Previous Work}\label{sec2:previous_work}

Research in the area of powered descent guidance began in earnest when the Apollo program set its sights on achieving a manned lunar landing. Until that point, no rocket-powered spacecraft had achieved a controlled soft landing on a celestial body. While the Apollo lunar modules were piloted by humans, the guidance system was capable of landing the spacecraft autonomously \cite{Klumpp1971}. Apollo guidance is a 3-degree-of-freedom (3-DoF) translation guidance method, with attitude commands computed separately using the generated thrust commands. Despite the popularity of powered descent guidance as a modern research topic, Apollo-derived methods still form the core of many guidance controllers~\cite{Mendeck2011,Sostaric2005,Singh2007}. 

While Apollo guidance is not, strictly speaking, optimal in any sense, its designers were aware of theoretical work that studied the problem of fuel-optimal soft landings~\cite{Klumpp1971,Lawden1963,Meditch1964,Edelbaum1969}. The work of Lawden in~\cite{Lawden1963} is often credited as the originator of several key insights, as he considered the general guidance problem and developed analytical expressions for optimal trajectories. At the time, however, numerical solutions to such problems were difficult to obtain. The first tractable solutions adequate for a flight implementation were limited to one-degree-of-freedom (1-DoF) vertical descent trajectories~\cite{Meditch1964}. These early results using the calculus of variations and Pontryagin's maximum principle were quite promising and continued to be developed after the Apollo program~\cite{Kornhauser1972,Breakwell1975,Marec1979,Azizov1986}.

Research on optimal powered descent garnered renewed interest due to unmanned Martian landing missions around the turn of the century. In~\cite{DSouza1997}, a simplified (no mass variable) 3-DoF landing problem was investigated using minimum acceleration and final time as the cost. No path constraints were imposed on state or control variables, but an analytic feedback control law was realized. Refs.~\cite{Topcu2005,Najson2005,Topcu2007} presented results for the 3-DoF landing problem that largely mirror Lawden's, though numerical trials using newly available software were provided to confirm theoretical results. Several authors have continued to look for analytical solutions to the first-order necessary conditions for the 3-DoF problem~\cite{Rea2010,Lu2018}. Very efficient numerical routines for obtaining optimal thrust programs have been developed~\cite{Lu2012,Lu2018}, though generally at the expense of state constraints.

At the same time, Refs.~\cite{Acikmese2005,Ploen2006,Acikmese2007a} provided an alternate viewpoint on the 3-DoF problem. These works took a convex programming approach, and showed that a non-convex lower bound on thrust magnitude can be relaxed to a convex constraint by using a slack variable. Using Pontryagin's maximum principle, they showed that in fact this relaxation is lossless, and the relaxed problem yields the same optimal solution as the original problem. A subsequent change of variables and an approximation led to a fully convex problem formulation that can be solved efficiently. This lossless convexification result bridged the remaining gap between theoretical understanding of the 3-DoF guidance problem and the ability to obtain numerical solutions quickly and efficiently. The theory of lossless convexification has since been expanded to non-convex thrust pointing constraints~\cite{Carson2011a,Carson2011b,Acikmese2013}, minimum-landing error problems~\cite{Blackmore2010a} and more general optimal control problems~\cite{Acikmese2011,Blackmore2012b,Harris2014}. These efforts culminated in the development of a custom second-order cone programming solver~\cite{Dueri2014,Dueri2016} and successful tests of the guidance system in terrestrial flight tests~\cite{Scharf2014,Scharf2017}.

More recent work has explored the generalized 6-DoF landing problem that considers both translational and rotational motion of a rigid body. These include the use of Lyapunov techniques~\cite{Lee2012c}, model predictive control~\cite{Lee2015a,Lee2017} and more recently feedforward trajectory generation techniques~\cite{Szmuk2016,Szmuk2017,Szmuk2018,SzmukReynolds2018,Reynolds2019,Szmuk2019}. The state variables selected for the 6-DoF problem formulation can be used to classify various methods. For example, one may use standard Cartesian variables in conjunction with unit quaternions, homogeneous transformations, or dual quaternions. We view this distinction as a design choice, but note that dual quaternions are an elegant and efficient parameterization for cases with constraints that couple rotation and translation. While a complete characterization of fuel optimal solution(s) for the 6-DoF problem is an active area of research, numerical techniques have yielded insightful results. All reported 6-DoF guidance schemes devise iterative strategies to obtain feasible solutions to nonlinear and non-convex optimal control problems, while approximating local optimality.  However, current research offers promising results that locally optimal solutions can be found by sequentially solving convex optimization problems with guaranteed convergence properties~\cite{Mao2016a,Mao2017,Mao2018,Virgili-Llop2018}. 

\subsection{Contributions}\label{sec2:contributions}

The contributions of this paper are: (i) a complete exposition of the dual quaternion-based 6-DoF powered descent guidance problem, (ii) the use of newly defined \textit{state-triggered constraints} to obtain trajectories that respect line of sight constraints during specific segments of a descent trajectory, (iii) the presentation of a real-time implementable algorithm and Monte Carlo study that examines its viability for on-board use. This paper provides open-loop error analysis and discusses how key algorithm parameters relate to open-loop accuracy. Runtime analysis on a 3.2 GHz Intel Core i5 processor is presented for the entire algorithm to provide an accurate estimate of computational capabilities.

Previous works such as~\cite{Lee2012c,Lee2017,Filipe2015,Reynolds2018} have proposed unit dual quaternion approaches using various feedback control techniques. In contrast, to the best of our knowledge, this work is the first to present a unit dual quaternion approach to feedforward trajectory generation. State-triggered constraints are used in this work to model range-triggered line of sight constraints that are present during Hazard Detection and Avoidance (HDA) or safe landing site selection. During HDA, optical sensors must be pointed towards certain locations of the ground in order to provide the navigation system with sufficient information to select a safe landing site. The guidance and control system must ensure that this constraint is sufficiently met. However, once a safe landing site is chosen, there may no longer be a need to point the optical sensor directly at it. As such, we view this as a constraint that should be disabled after either a fixed amount of time or below a certain range from the landing site. As we show in~\sref{sec3:computation_performance}, this work is the first to fuse such a requirement into a guidance method that is conducive to real-time implementation. 

This paper is organized as follows. First,~\sref{sec:dqs} introduces the quaternion and dual quaternion formalism used in this work, while~\sref{sec2:rigid_body_motion} specializes these to 6-DoF rigid body motion. Next,~\sref{sec:problem_statement} details the non-convex free-final time optimal control problem, while~\sref{sec:sol_method} describes the algorithm developed to solve the aforementioned problem. Two numerical case studies are presented in~\sref{sec:numerical_demo} that highlight the primary contributions of this work and capabilities of the algorithm. Lastly,~\sref{sec:conclusion} offers concluding remarks and directions for future research.

\section{Dual Quaternions and Rigid Body Motion}\label{sec:dqs}

The set of possible orientations of one coordinate frame relative to another is given by the special orthogonal group in three dimensions, $SO(3) \definedas \{ C\in\real^{3\times 3} \,|\, \det{C}=+1,\,C^{-1}=C\tran\}$. The elements of this set are commonly referred to as \textit{rotation matrices} and have six free parameters. On the other hand, a quaternion $\q\in\setq$ has four parameters and is composed of a three parameter vector part, $\qv$, and a scalar part, $\qs$, such that $\q=(\qv,\qs)$. The set of unit quaternions is denoted by $\setuq \definedas \{ \q\in\setq \,|\, \q\ipQ\q = 1\}$, where $\ipQ$ denotes the Euclidean inner product in~$\setq$. Unit quaternions provide a minimal singularity-free representation of the set of rotation matrices and form a double cover of $SO(3)$. The set $\setuq$ can be thought of as the three-sphere embedded within the quaternion manifold $\setq$.

Dual quaternions extend the quaternion parameterization to capture both the relative orientation and relative position of two coordinate frames. Dual quaternions can be elegantly derived using Clifford algebras as in~\cite{McCarthy1990}, or by geometric construction as in the original work~\cite{Clifford1882}.
We denote a dual quaternion by
\begin{equation}
\dqq = \q_1 + \eps \q_2 \, \in \setdq\,,
\label{eq:dq_setdefn}
\end{equation} 
where $\q_1,\,\q_2 \in \setq$ are quaternions, and $\eps\ne0$ is termed the \textit{dual unit} that satisfies the property $\eps^2 = 0$. We call $\q_1$ the \textit{real part} and $\q_2$ the \textit{dual part} of the dual quaternion $\dqq$. Under the dual quaternion dot product~\cite{McCarthy1990}, a unit dual quaternion satisfies
\begin{equation}
\dqq \ipDQ \dqq = \left( \q_1 + \eps \q_2 \right) \ipDQ \left( \q_1 + \eps \q_2 \right) = \q_1 \ipQ \q_1 + \eps \left( \q_1 \ipQ \q_2 + \q_2 \ipQ \q_1 \right) = 1 + \eps 0.
\label{eq:dq_dot_dq}
\end{equation}
It follows then that the real and dual parts of $\dqq$ must satisfy $\q_1 \ipQ \q_1 = 1$ and $\q_1 \ipQ \q_2 = 0$. Consequently, we define the set of unit dual quaternions as
\begin{equation}
\setudq \definedas \left\lbrace \dqq = \q_1 + \eps \q_2 \,|\, \q_1 \ipQ \q_1 = 1,\; \q_1 \ipQ \q_2 = 0 \right\rbrace.
\label{eq:udq_setdefn}
\end{equation}

Unit dual quaternions form a submanifold within the dual quaternion manifold. Note that the first constraint forces the real part of a unit dual quaternion to be an element of the three-sphere $\setuq$. The second constraint forces the dual part of a unit dual quaternion to be an element of the (three-dimensional) tangent plane of the three-sphere at the point $\q_1$. As such, the set of unit dual quaternions in~\eqref{eq:udq_setdefn} can be interpreted as the union of the Cartesian products of each $\q_1\in\setuq$ and the corresponding tangent plane to the three-sphere at $\q_1$. 

\subsection{Dual Quaternion Operations}\label{sec2:dq_operations}
Let $\dq{a},\dq{b} \in \setudq$ be two unit dual quaternions, and let $\quat{a},\quat{b} \in \setuq$ be two unit quaternions. We define quaternion multiplication using Hamilton's convention as
\begin{equation}
\quat{a} \otimes \quat{b} = \left( a_4 \quat{b}_v + b_4 \quat{a}_v + \quat{a}_v^{\times} \quat{b}_v ,\; a_4 b_4 - \quat{a}_v \ipQ \quat{b}_v \right),
\label{eq:q_mult}
\end{equation}
and dual quaternion multiplication as
\begin{equation}
\dq{a} \otimes \dq{b} = \quat{a}_1 \otimes \quat{b}_1 + \eps \left( \quat{a}_1 \otimes \quat{b}_2 + \quat{a}_2 \otimes \quat{b}_1 \right).
\label{eq:dq_mult}
\end{equation} 

Note that the product $\dq{a}\otimes\dq{b}$ is also a unit dual quaternion. Next, the quaternion cross product is defined as
\begin{equation}
\quat{a} \oslash \quat{b} = \left( a_4 \quat{b}_v + b_4 \quat{a}_v + \quat{a}_v^{\times} \quat{b}_v,\; 0 \right),
\label{eq:q_cross}
\end{equation}
which is used in turn to define the dual quaternion cross product~\cite{Lee2017,Lee2015a}
\begin{equation}
\dq{a} \oslash \dq{b} = \quat{a}_1 \oslash \quat{b}_1 + \eps \left( \quat{a}_1 \oslash \quat{b}_2 + \quat{a}_2 \oslash \quat{b}_1 \right).
\label{eq:dq_cross}
\end{equation}
We define the quaternion conjugate as $\quat{a}^* = (-\quat{a}_v,\,a_4)$, and the dual quaternion conjugate as
\begin{equation}
\dq{a}^* = \quat{a}_1^* + \eps \quat{a}_2^*.
\label{eq:dq_conj}
\end{equation}

To use more familiar matrix-vector analysis, we embed the set of unit dual quaternions in the eight-dimensional Euclidean space $\real^8$. Denoting the Euclidean inner product in $\real^8$ using $\cdot\tran\cdot\,$, and using the natural isomorphism
\begin{equation}
\dq{a} = \quat{a}_1 + \eps \quat{a}_2 \in \setudq \quad \mapsto \quad \dq{a} = \begin{bmatrix}
\quat{a}_1 \\ \quat{a}_2
\end{bmatrix} \in \setudqR \definedas \left\lbrace \dq{a} \in \real^8 \,|\, \quat{a}_1\tran \quat{a}_1 = 1 \ \text{and} \ \quat{a}_1\tran \quat{a}_2 = 0 \right\rbrace,
\label{eq:dq_2_8vec}
\end{equation}
we henceforth view the set of unit dual quaternions as the six-dimensional subset $\setudqR \subset \real^8$. By virtue of the constraint imposed on the first four elements of $\dq{a} \in \setudqR$, we view unit quaternions as elements of the three-dimensional subset $\setuqR\subset\real^4$. As a result, we represent the operation in~\eqref{eq:q_mult} with the following matrix expressions
\begin{equation}
\quat{a} \otimes \quat{b} = \qskew{\quat{a}}\,\quat{b} = \qskewstar{\quat{b}}\,\quat{a}
\label{eq:q_mult_mat}
\end{equation}
where,
\begin{equation*}
\qskew{\quat{a}} \definedas \begin{bmatrix}
a_4 \eye{3} + \quat{a}_v^{\times} & \quat{a}_v \\ -\quat{a}_v^T & a_4
\end{bmatrix} \quad \text{and} \quad \qskewstar{\quat{b}} \definedas \begin{bmatrix}
b_4 \eye{3} -\quat{b}_v^{\times} & \quat{b}_v \\ - \quat{b}_v^T & b_4
\end{bmatrix}.
\end{equation*}
Using these definitions, we can then rewrite~\eqref{eq:dq_mult} as
\begin{equation}
\dq{a} \otimes \dq{b} = \qskew{\dq{a}}\,\dq{b} = \qskewstar{\dq{b}}\dq{a},
\label{eq:dq_mult_mat}
\end{equation}
where,
\begin{equation*}
\qskew{\dq{a}} \definedas \begin{bmatrix}
\qskew{\quat{a}_1} & \zeros{4}{4} \\ \qskew{\quat{a}_2} & \qskew{\quat{a}_1}
\end{bmatrix} \quad \text{and} \quad \qskewstar{\dq{b}} \definedas \begin{bmatrix}
\qskewstar{\quat{b}_1} & \zeros{4}{4} \\ \qskewstar{\quat{b}_2} & \qskewstar{\quat{b}_1}
\end{bmatrix}.
\end{equation*}

The matrices in~\eqref{eq:dq_mult_mat} are structured so that the matrix-vector multiplication gives the same result as the definition in~\eqref{eq:dq_mult}. As a result, the dual unit is no longer present in these expressions. The columns of these matrices may also be interpreted as the (non-commutative) projection of the dual quaternion onto the basis vectors of the Clifford sub-algebra used to derive them~\cite{McCarthy1990}.  The quaternion and dual quaternion cross products can be rewritten as matrix-vector products using the same methods; see~\cite{Lee2015a,Lee2017} for more details.  



\subsection{Rigid Body Motion}\label{sec2:rigid_body_motion}

\begin{figure}[t!]
\centering

\begin{tikzpicture}[scale=1]
    \isoaxesNL{0}{0}{90}{1}{1}{$\Finertial$}
    
    \isoaxesNLC{0}{0}{90+25}{1}{1}{}{gray}
    
    \coordinate (rx) at ( 0.85,0.25,0 );
    \coordinate (ry) at ( -0.25,0.85,0);
    \coordinate (rz) at ( -0.05,-0.25,0.6 );
    \coordinate (c) at (2,-2,0);
    
	\isoaxesNLC{2}{-2}{90+25}{1}{1}{}{gray}
    \node[gray] at (2,-2.5) {$\Fbody$};
    
    \draw[gray,thick,domain=-30:120,dashed,->] plot ({0.5*cos(\x)}, {0.5*sin(\x)}) node [right, xshift=0.4cm, yshift=0.1cm] {$\bm{q}$};
    \draw[thick,dashed,->] (0,0,0) -- (c) node[pos=0.5,above right] {$\rB$};
    
    \node[thick] at (3,0,0) {$\equiv$};
    \node[] at (3,-0.30,0) {$\dqq$};
    
    \begin{scope}[shift={(5,0,0)}]
        	\isoaxesNL{0}{0}{90}{1}{1}{$\Finertial$}
    		\coordinate (c) at (2,-2,0);    
    		\isoaxesNL{2}{-2}{90}{1}{1}{}  
    		\isoaxesNLC{2}{-2}{90+25}{1}{1}{}{gray} 
    		\node[gray] at (2,-2.5) {$\Fbody$};
    		
   		\draw[gray,thick,domain=-30:120,dashed,->] plot ({2+0.5*cos(\x)}, {-2+0.5*sin(\x)}) node [right, xshift=0.4cm, yshift=0.1cm] {$\bm{q}$};
   		
    		\draw[thick,dashed,->] (0,0,0) -- (c) node[pos=0.5,above right] {$\rI$};
    \end{scope}
    
\end{tikzpicture}

\caption{Rotation followed by translation is equivalent to a translation followed by a rotation.}
\label{fig:frame-diff}
\end{figure}
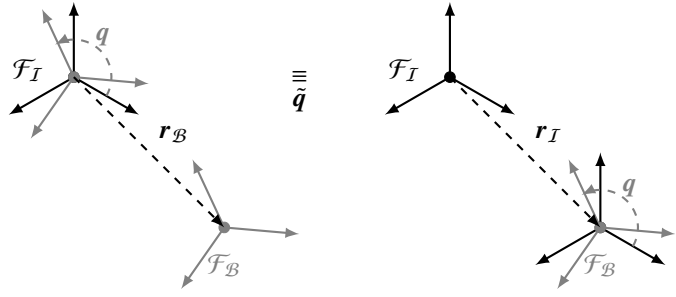

Let us consider two three-dimensional coordinate frames: $\Finertial$ an inertial frame with its origin at the landing site, and $\Fbody$ a body-fixed frame whose origin is the center of mass of the vehicle. Pure rotation by a unit quaternion $\q$ is represented by the dual quaternion $\q + \eps \zeros{4}{1}$. Pure translation may be described either by $\qid + \eps \frac{1}{2}\rI$ or $\qid + \eps \frac{1}{2}\rB$. Note that $\rI$ and $\rB$ are assumed to have a zero scalar part and treated as \textit{pure} quaternions here. To map between the inertial and body coordinate frames, we may either perform a translation by $\rI$ followed by a rotation by $\q$, or perform a rotation by $\q$ followed by a translation $\rB$. These two sequences are geometrically equivalent, as shown in Fig.~\ref{fig:frame-diff}.

The composition of a rotation and translation is represented by dual quaternion multiplication. Using the definitions of pure rotation and translation, and taking~\eqref{eq:dq_2_8vec} into account, it follows that the unit dual quaternion representing the difference between these two frames is 
\begin{equation}
\dqq = \begin{bmatrix}
\q \\ \frac{1}{2} \rI \otimes \q
\end{bmatrix} = \begin{bmatrix}
\q \\ \frac{1}{2} \q \otimes \rB
\end{bmatrix} \in \setudqR.
\label{eq:dq_rigidbody_defn}
\end{equation}
The first expression describes a translation by $\rI$ followed by a rotation $\q$, whereas the second expression describes a rotation by $\q$ followed by a translation $\rB$. 
The equivalence of the two expressions in~\eqref{eq:dq_rigidbody_defn} leads to the observation that
\begin{equation}
\rI = \q \otimes \rB \otimes \q^* \quad \text{and} \quad \rB = \q^* \otimes \rI \otimes \q.
\label{eq:i2b_quat}
\end{equation}

Velocities can be represented using dual quaternions by appending a zero scalar part to the usual angular and linear velocities such that $\wB,\vB \in \real^4$. The \textit{dual velocity} is then defined to be
\begin{equation}
\dqw = \wB + \epsilon \vB \in \setdq \quad \mapsto \quad
\dqw = \begin{bmatrix}
\wB \\ \vB 
\end{bmatrix} \in \real^8.
\label{eq:dqw_defn}
\end{equation}
%

Our definition of the dual quaternion in~\eqref{eq:dq_rigidbody_defn} in conjunction with the quaternion triple identity~\cite{Lee2015a} leads to the following lemma, which can be proved by construction.
\begin{lemma}\label{prop:quat_id_quadratic}
Let $\q\in\setuqR$ be the unit quaternion that maps inertial coordinates to body coordinates. Let $\bm{r}_{\inertial}$ and $\bm{r}_{\body}$ be the coordinates of the position vector in $\Finertial$ and $\Fbody$ respectively. Consider the pure quaternions $\bm{a},\,\bm{b},\in\real^4$ with coordinates $\bm{a}_{\inertial},\,\bm{b}_{\inertial}$ in the inertial frame and $\bm{a}_{\body},\,\bm{b}_{\body}$ in the body frame.  Using~\eqref{eq:dq_rigidbody_defn}, the following equations hold:
\begin{subequations}
\begin{alignat}{3}
        \bm{r}_{\inertial}\tran \bm{a}_{\inertial} &= \dqq\tran M_1 \dqq, \quad &\text{where} \quad  M_1 &= \begin{bmatrix}
        \zeros{4}{4} & \qskew{\bm{a}_{\inertial}}\tran \\ \qskew{\bm{a}_{\inertial}} & \zeros{4}{4}
        \end{bmatrix}, \label{eq:quat_id_II} \\
    \bm{r}_{\body}\tran \bm{a}_{\body} &= \dqq\tran M_2 \dqq, \quad  &\text{where} \quad M_2 &= \begin{bmatrix}
        \zeros{4}{4} & {\qskewstar{\bm{a}_{\body}}}\tran \\ \qskewstar{\bm{a}_{\body}} & \zeros{4}{4}
        \end{bmatrix}, \label{eq:quat_id_BB}\\
    \bm{a}_{\inertial}\tran \bm{b}_{\inertial} &= \dqq\tran M_3 \dqq, \quad  &\text{where} \quad M_3 &= \begin{bmatrix}
        \qskew{\bm{a}_{\inertial}}\qskewstar{\bm{b}_{\body}} & \zeros{4}{4} \\ \zeros{4}{4} & \zeros{4}{4}
        \end{bmatrix}, \label{eq:quat_id_IB}  \\ 
    \| \rI \|_2 &= \| 2 E_d \dqq \|_2,  \quad &\text{where} \quad  E_d &= \begin{bmatrix}
    \zeros{4}{4} & \zeros{4}{4} \\ \zeros{4}{4} & \eye{4}
    \end{bmatrix}. \label{eq:norma_equiv_dq}
\end{alignat}
\end{subequations}
\end{lemma}

\subsubsection{Kinematics \& Dynamics}\label{sec3:kinematics_dynamics}

We assume in this work that the vehicle's mass varies as a function of thrust, but we neglect changes in the inertia matrix and center of mass. In previous work, we investigated the effects of variable inertia on the trajectories generated using similar methods to those described here~\cite{Reynolds2019}. In the scenarios studied in~\cite{Reynolds2019}, trajectories did not deviate significantly from those obtained with a constant inertia matrix, and so these variations are not considered here. The mass, $m \in \realpp$, is assumed to vary as a linear function of the thrust magnitude according to
\begin{equation}
\dot{m} = - \alpha \| \uB \|_2, \qquad \alpha \definedas \frac{1}{I_{\text{sp}} g_e},
\label{eq:mass_dyn}
\end{equation}
where $\uB \in \real^3$ is the thrust vector in the body frame.

The dual quaternion kinematic equation can be obtained by directly computing the time derivative of \eqref{eq:dq_rigidbody_defn}. Using the well-known quaternion kinematic equation $\dot{\q}=\frac{1}{2}\q \otimes \wB$, we write
\begin{align}
\der{\dqq}{t} = \frac{d}{dt} \begin{bmatrix}
\q \\ \frac{1}{2} \rI \otimes \q
\end{bmatrix}  = \begin{bmatrix}
\qdot \\ \frac{1}{2} ( \rIdot \otimes \q + \rI \otimes \qdot )
\end{bmatrix} &= \begin{bmatrix}
\frac{1}{2} \q \otimes \wB \\ \frac{1}{2} ( \q \otimes \vB + \frac{1}{2} \rI \otimes \q \otimes \wB )
\end{bmatrix} = \frac{1}{2} \dqq \otimes \dqw 
\label{eq:dq_kinematics}
\end{align}
where we have used~\eqref{eq:i2b_quat} to write $\vI \otimes \q = \q \otimes \vB$ and~\eqref{eq:dq_mult_mat} to obtain the last equality.

The translational control capabilities of landing vehicles are typically dominated by the rocket engine(s), whereas the attitude control capabilities are dictated by both the rocket engine(s) and/or a reaction control system (RCS). For guidance trajectory design, we assume that the rocket engine(s) provide both the translation and attitude control authority.  We assume that an RCS system is used strictly for closed-loop attitude control. However, we stress that RCS can be easily added to the problem formulation.

The dynamics are obtained by using the Newton-Euler equations in a rotating frame. Taking a derivative of the linear and angular momenta in the rotating body frame leads to
\begin{subequations}\label{eq:newton_euler_eqms}
\begin{align}
\der{}{t}\big(m \vB \big) &= m \vBdot + \wB^{\times} ( m \vB ) = \sum \FB + \uB, \label{eq:newton_eqm} \\
\der{}{t}\big(J \wB \big) &= J \wBdot + \wB^{\times} ( J \wB ) = \sum \TB + \ru^{\times} \uB ,  \label{eq:euler_eqm}
\end{align}
\end{subequations}
where vector $\ru\in\real^3$ denotes the constant body-frame vector from the vehicle's center of mass to the point where the thrust is applied by a single rocket engine\footnote{Multiple engines can easily be incorporated in our framework by appropriate redefinition of the $\uB$ terms in~\eqref{eq:newton_euler_eqms}.}. We assume that $\sum \FB = m \gB$ and $\sum \TB = 0$, where $m\gB$ is the force due to gravity in the body frame. Note that momentum changes due to mass variability are accounted for in~\eqref{eq:newton_euler_eqms} by our definition of the mass depletion dynamics in~\eqref{eq:mass_dyn} and $\uB$; see~\cite{Thomson1986} for details.

Combining the expressions in~\eqref{eq:newton_euler_eqms} with the definition of the dual velocity in~\eqref{eq:dqw_defn} allows us to express the dual quaternion dynamics as
\begin{equation}
\J \dqwdot + \dqw \oslash \J \dqw = \Phi \uB + m\dqgB,
\label{eq:dq_dynamics}
\end{equation}
where,
\begin{equation*}
\bm{J} \definedas \left[ \begin{array}{c|c}
0_{4 \times 4} & \begin{array}{cc}
m I_3 & \zeros{3}{1} \\ \zeros{1}{3} & 1
\end{array} \\
\hline
\begin{array}{cc}
J & \zeros{3}{1} \\ \zeros{1}{3} & 1
\end{array} & 0 _{4 \times 4 }
\end{array} \right]_{8 \times 8} \qquad \Phi \definedas \left[ \begin{array}{cc}
I_3 \\ \zeros{1}{3} \\
\hline
\ru^{\times} \\ \zeros{1}{3}
\end{array} \right]_{8\times 3} \qquad \dq{g}_{\body} \definedas \left[ \begin{array}{c}
\gB \\ \zeros{4}{1}
\end{array} \right]_{8\times 1}.
\end{equation*}
Note that the \textit{dual inertia} matrix $\bm{J}$ is always invertible. We refer the reader to~\cite{Lee2017,Lee2015a,Filipe2015,Reynolds2018} for more details on rigid body kinematics and dynamics using dual quaternions.

\section{Problem Statement}\label{sec:problem_statement}

This section formulates the non-convex free-final-time powered descent guidance problem that is considered in this paper. Having already stated the equations of motion, this section focuses specifically on the state and control constraints imposed during powered descent. Without loss of generality, we assume that the inertial frame has its origin at the nominal landing site and is constructed from the orthonormal vectors $\{\xI,\,\yI,\,\zI\}$. These vectors are oriented such that $\xI$ represents the downrange direction, $\yI$ represents the crossrange direction, and $\zI$ points locally up. Similarly, the body frame $\Fbody$ has its origin at the vehicle's (constant) center of mass and is constructed from the orthonormal vectors $\{\xB,\,\yB,\,\zB\}$, where $\zB$ is chosen to point along the vehicle's vertical axis. We assume that these vectors coincide with the vehicle's principal axes of inertia.

This section is organized as follows. A \textit{baseline} set of state and control constraints and boundary conditions are first described in~\sref{sec2:baseline}. Next, state-triggered constraints are introduced in~\sref{sec2:stc} and their application to slant-range-triggered line of sight constraints is presented. Lastly, a complete statement of the powered descent guidance problem is given in~\sref{sec2:ncvx_problem_statement}.

\subsection{Baseline Problem Constraints}\label{sec2:baseline}



\subsubsection{State Constraints}\label{sec3:baseline_state}

To ensure that trajectories do not use more fuel than is stored on-board, a constraint is enforced on the mass of the vehicle according to
\begin{equation}
m \geq \mdry.
\label{eq:max_fuel}
\end{equation}


Next, we define the approach cone angle to be the angle formed between $\rI$ and $\zI$. An approach cone constraint is used to ensure the vehicle's position lies above the surface of the planet, while also ensuring sufficient elevation at large distances from the landing site. This constraint can be expressed as
%
\begin{equation}
-\rI\tran \zI +  \| \rI \|_2 \cos \gsmax \leq 0,
\label{eq:gs_max}
\end{equation}
where we assume that $\gsmax\in[0\dg,90\dg]$.

The approach cone constraint~\eqref{eq:gs_max} is expressed in terms of the dual quaternion by using~\eqref{eq:quat_id_II} and~\eqref{eq:norma_equiv_dq} to write
\begin{equation}
\gsf(\dqq) \definedas -\dqq\tran \gsM \dqq + \| 2 E_d \dqq \|_2 \cos \gsmax \leq 0.
\label{eq:dq_gs_max}
\end{equation}
It was shown in~\cite{Lee2015a,Lee2017} that $\gsf:\setudqR \rightarrow \real$ is convex over the bounded domain $\dom{\gsf} = \big\{ \dqq\in\setudqR \,|\, \dqq\tran \dqq \leq 1 + \frac{1}{4}\Delta^2 \big\}$, where $\| \rI\|_2 \leq \Delta$ is an upper bound on the distance from the landing site.

The vehicle is also subject to a tilt angle constraint that limits the angle formed between the vehicle's vertical axis $\zB$ and the inertial direction, $\zI$. When both vectors are viewed inertially, we can formulate the tilt angle constraint as
\begin{equation}
    -\zI\tran \left( \q \otimes \zB \otimes \q^* \right) + \cos \tiltmax \leq 0.
    \label{eq:tilt_max}
\end{equation}
The constraint~\eqref{eq:tilt_max} is expressed in terms of the dual quaternion by using~\eqref{eq:quat_id_IB} as
\begin{equation}
    \tiltf(\dqq) \definedas \dqq\tran  \tiltM \dqq + \cos \tiltmax \leq 0.
    \label{eq:dq_tilt_max}
\end{equation}
It was shown in~\cite{Lee2017,Reynolds2019} that for any $\tiltmax\in[0,90]\dg$, the function $\tiltf$ is convex for all $\dqq\in\setudqR$.
The last state constraint that we consider bounds the allowable angular rates by  enforcing
\begin{equation}
    \| \wB \|_{\infty} \leq \wmax.
    \label{eq:w_max}
\end{equation}

\subsubsection{Control Constraints}\label{sec3:baseline_ctrl}

We have assumed in~\eqref{eq:dq_dynamics} that a single rocket engine provides both the translation and attitude control authority. At times, rocket engines necessitate operation in restricted thrust and gimbal angle regimes. For example, the main engines of Apollo-era landers could either operate at $93\%$ thrust or in the permitted interval of $11\%$ to $65\%$ of the rated thrust value~\cite{Klumpp1971}. The forbidden thrust regions were avoided to prevent cavitation in the propulsion system. To model permitted thrust regions, we place restrictions on the norm of the thrust vector according to
\begin{equation}
\umin \leq \| \uB \|_2 \leq \umax,
\label{eq:thrust_bound}
\end{equation}
where $[\umin,\umax]\subset\realpp$ denotes the permitted thrust interval. 

We assume the engine may be gimbaled symmetrically about two axes and define the gimbal angle to be the total angular deflection of the thrust vector from its nominal position. We model the mechanical limitations of the engine using the gimbal angle constraint given by
%
\begin{equation}
\| \uB \|_2 \leq \sec \gimbalmax \bm{z}_{\body}\tran \uB,
\label{eq:max_gimbal}
\end{equation}
where we assume that $\gimbalmax \ll 90\dg$. 

The last constraint imposed on the control is a rate constraint that ensures commanded thrust vectors do not change too rapidly for the engine to follow. This can be formulated as
\begin{subequations} \label{eq:max_du}
\begin{gather} 
\|E_{xy}\duB\|_2 \leq \dgimmax \zB\tran \uB, \label{eq:max_gimbal_rate} \\
-\dot{u}_{z,\max} \leq \zB\tran \duB \leq \dot{u}_{z,\max}, \label{eq:max_throttle_rate}
\end{gather}
\end{subequations}
where $E_{xy}\in\real^{2\times 3}$ selects the thrust vector components in the $\xB$-$\yB$ plane.

\subsubsection{Boundary Conditions}\label{sec3:baseline_bcs}

We assume that the initial mass, position, velocity and angular velocity are fixed. We choose to let the optimization select the initial attitude as part of the optimization for two reasons. First, typical lunar descent sequences have a short pitch-up attitude maneuver that occurs immediately prior to the final approach phase~\cite{Klumpp1971,Sostaric2005}. A guidance scheme like the one described herein would select the attitude to be achieved at the end of the pitch-up maneuver, thereby ensuring the initial attitude satisfies constraints important to the subsequent powered descent maneuver. Second, we have found empirically that leaving this variable free offers better convergence behaviour over a wider range of initial conditions, in particular when line of sight or pointing constraints are incorporated into the formulation (see~\sref{sec3:stc_los}). 
These constraints are enforced as equality constraints according to
\begin{equation}
    m(\ti) = \mi, \quad \dqq(\ti) = \dqqig(\q(\ti)) \definedas \begin{bmatrix} \q(\ti) \\ \frac{1}{2} \ri{\inertial}\otimes\q(\ti) \end{bmatrix} ,\quad \dqw(\ti) = \dqwig(\q(\ti)) \definedas \begin{bmatrix} \wi \\ \q^*(\ti) \otimes \vi{\inertial} \otimes \q(\ti) \end{bmatrix},
    \label{eq:bcs_init}
\end{equation}
where $\mi$, $\ri{\inertial}$, $\vi{\inertial}$ and $\wi$ represent the specified initial mass, inertial position, inertial velocity and angular velocity. 
%
At the final time, we fix the dual quaternion and dual velocity to prescribed points, but leave the final mass free. These are enforced as equality constraints according to
\begin{equation}
    \dqq(\tf) = \dqq_f, \quad \dqw(\tf) = \dqw_f,
    \label{eq:bcs_final}
\end{equation}
where $\dqq_f$ and $\dqw_f$ are constructed from the target inertial position, $\rf{\inertial}$, attitude $\qf$, inertial velocity, $\vf{\inertial}$, and angular velocity, $\wf$ using the definitions in~\sref{sec:dqs}. 

\subsection{State-Triggered Constraints}\label{sec2:stc}

State-triggered constraints were recently introduced in~\cite{SzmukReynolds2018,Szmuk2019,Reynolds2019}. Simply stated, they allow discrete decisions to be formulated in a continuous optimization framework (e.g., successive convexification). In what follows, we provide a brief overview of state-triggered constraints and then show how they are used to formulate a slant-range-triggered line of sight (LoS) constraint for the powered descent guidance problem. We highlight that the ``state'' referred to here is not limited to the traditional definition from linear systems theory. Rather, it connotes a (set of) variable(s) from an optimization problem and applies equally to control- or time-triggered constraints.

\subsubsection{Logical Statement}

Each state-triggered constraint is composed of a \textit{trigger condition} and a \textit{constraint condition}. The trigger condition is given by the inequality~$\stcg(\stcz)<0$, where~$\stcz\in\real^{\stcnz}$ denotes the solution variable of a parent optimization problem, and~$\stcg(\cdot):\real^{\stcnz}\rightarrow\real$ is called the \textit{trigger function}. The constraint condition is given by the inequality~$\stcc(\stcz)\leq 0$, where~$\stcc(\cdot):\real^{\stcnz}\rightarrow\real$ is called the \textit{constraint function}. We assume that both $\stcg(\cdot)$ and $\stcc(\cdot)$ are differentiable with respect to their arguments. A state-triggered constraint is defined as the following logical implication:
\begin{equation}
    \stcg(\stcz) < 0\;\Rightarrow\;\stcc(\stcz)\leq 0.
    \label{eq:stc_logic}
\end{equation}
The practical value of~\eqref{eq:stc_logic} is also evident from the contrapositive implication: the constraint condition may be violated \textit{only if} the trigger condition is not satisfied.

\subsubsection{Continuous Formulation}

Mixed-integer techniques are the most common way to implement constraints like~\eqref{eq:stc_logic}, and require the introduction of discrete decision variables. Despite the existence of efficient branch and bound algorithms, these approaches suffer from worst-case exponential computational complexity~\cite{Biegler2014,Richards2015}. This computational complexity is further compounded by the iterative nature of the solution process required to solve the problem addressed in this paper; each combination of discrete decisions in the branch and bound sequence would require multiple iterations to converge. Consequently, we argue that mixed-integer solution strategies are not well-suited for solving powered descent guidance problems in real-time.

In contrast, a continuous variable formulation of~\eqref{eq:stc_logic} can leverage the iterations of a sequential approach to embed binary decisions without resorting to discrete decision variables, thus avoiding a large computational penalty. A continuous variable formulation of~\eqref{eq:stc_logic} is given by
\begin{equation}
    \stch(\stcz) \definedas \shat(\stcz)\,\stcc(\stcz) \leq 0,
    \label{eq:stc_cont}
\end{equation}
where~$\shat(\stcz)\definedas -\min\big(0,\stcg(\stcz)\big)$. We refer the reader to~\cite{SzmukReynolds2018} for more details regarding~\eqref{eq:stc_cont}. Evidently, when the trigger condition is satisfied (i.e.,~$\stcg(\stcz) < 0$), then~$\shat(\stcz)>0$ and~\eqref{eq:stc_cont} implies that~$\stcc(\stcz)\leq 0$. When the trigger condition is not satisfied, then~$\shat(\stcz)=0$ and~\eqref{eq:stc_cont} is trivially satisfied for any value of~$\stcc(\stcz)$. Thus, we conclude that~\eqref{eq:stc_logic} and~\eqref{eq:stc_cont} are logically equivalent.\footnote{The state-triggered constraint in~\eqref{eq:stc_logic} and~\eqref{eq:stc_cont} can be reformulated by using an \textit{equality} constraint condition~$\stcc(\stcz)=0$, as in~\cite{SzmukReynolds2018}.}

\subsubsection{Compound State-Triggered Constraints}

Certain applications may have one constraint condition that is enforced when a combination of trigger conditions are satisfied or, conversely, have multiple constraint conditions that are enforced when one trigger condition is satisfied. Such constraints are called \textit{compound state-triggered constraints} and were introduced in~\cite{Szmuk2019}. The trigger and constraint conditions of compound state-triggered constraints are composed by using Boolean \textit{and} and \textit{or} operations. In this work we consider a compound-\textit{and} trigger condition and a single constraint condition. This constraint is logically expressed as
\begin{equation}
    \bigwedge_{i=1}^{\stcng}\big(\stcg_i(\stcz) < 0\big) \;\Rightarrow\; \stcc(\stcz)\leq 0,
    \label{eq:cstc_logic}
\end{equation}
where each~$\stcg_i(\cdot)$ is defined as in the scalar case. The continuous formulation corresponding to~\eqref{eq:cstc_logic} is given by
\begin{equation}
    \stch_{\land}(\stcz) \definedas \left[\prod_{i=1}^{\stcng}\shat_i(\stcz) \right]\cdot\stcc(\stcz) \leq 0,
    \label{eq:cstc_cont}
\end{equation}
where each~$\shat_i(\cdot)\definedas -\min(0,\stcg(\stcz))$ is defined as in the scalar case. On the one hand, if \textit{all} of the trigger conditions \textit{are} satisfied, then~$\prod_i\shat_i(\cdot) > 0$ and~\eqref{eq:cstc_cont} implies that~$\stcc(\stcz)\leq 0$. On the other hand, if \textit{any} of the trigger conditions \textit{are not} satisfied, then~$\prod_i\shat_i(\cdot) = 0$ and~\eqref{eq:cstc_cont} is trivially satisfied for any value of~$\stcc(\stcz)$. Hence, we conclude that~\eqref{eq:cstc_logic} and~\eqref{eq:cstc_cont} are logically equivalent.

\subsubsection{Application to Powered Descent: Slant-Range-Triggered Line of Sight}\label{sec3:stc_los}

\begin{figure}[b!]
    \centering
    \begin{tikzpicture}
    \tikzmath{\cxi=-5.5; \cyi=-2.5; 
              \cx1=\cxi+11; \cy1=\cyi+3.5; \aa1=-45; 
              \cx2=\cxi+6;  \cy2=\cyi+4.5; \aa2=-78; 
              \cx3=\cxi+1;  \cy3=\cyi+3.0; \aa3=-58; 
              \lrocket=1.6; \lflame=0.6;
              \Lcone=1.75; \Acone=20; \Dcone=0.2;
              \LL=3.0; \fact=1.4; \pfact=0.85; \tfact=0.65;
              \aac=20; 
              \aaac1=\aa1+\aac; 
              \aaac2=\aa2+\aac; 
              \aaac3=\aa3+\aac; 
              \aad=atan((\cy2-\cyi)/(\cx2-\cxi)); 
              \da=10; 
              \RRR1 = 4.40;
              \RRR2 = 9.5;
              \aaa1 = 15;
              \aaa2 = 65;
              \aaa3 = 40;
              \aaa4 = 15;
              \ppx1=\cxi+\RRR1*cos(\aaa1); \ppy1=\cyi+\RRR1*sin(\aaa1);
              \ppx2=\cxi+\RRR1*cos(\aaa2); \ppy2=\cyi+\RRR1*sin(\aaa2);
              \ppx3=\cxi+\RRR2*cos(\aaa3); \ppy3=\cyi+\RRR2*sin(\aaa3);
              \ppx4=\cxi+\RRR2*cos(\aaa4); \ppy4=\cyi+\RRR2*sin(\aaa4);
              \dang1=0; \aaaa1=5;
              \dang2=0; \aaaa2=11.5;
    }
    
    \definecolor{color0}{rgb}{0.00,0.50,0.00} 
    \definecolor{color1}{rgb}{0.00,0.65,0.00} 
    \definecolor{color2}{rgb}{1.00,0.00,0.00} 
    \definecolor{color3}{rgb}{0.85,0.00,1.00} 
    \definecolor{color4}{rgb}{0.00,0.00,1.00} 
    
    \pane{0}{0}{0}{16.5}{8.0}{0.5}{color=black,fill=beige!20}
    
    \fill[color3!10,opacity=1.0]
        (\ppx1,\ppy1) arc(\aaa1:\aaa2:\RRR1) --
        (\ppx2,\ppy2) to[out=90,in=170] (\ppx3,\ppy3) --
        (\ppx3,\ppy3) arc(\aaa3:\aaa4:\RRR2) --
        (\ppx4,\ppy4) to[out=270,in=-10] (\ppx1,\ppy1) --
        cycle;
    \draw[black,line width=0.5mm] (\ppx1,\ppy1) arc(\aaa1:\aaa2:\RRR1);
    \draw[black,line width=0.5mm] (\ppx3,\ppy3) arc(\aaa3:\aaa4:\RRR2);
    \draw[black,thin,dash pattern=on 0.025cm off 0.05cm]
        ({\cxi+\RRR1*cos(\aaa1-\dang1)},
         {\cyi+\RRR1*sin(\aaa1-\dang1)}) arc(\aaa1-\dang1:\aaa1-\dang1-\aaaa1:\RRR1);
    \draw[black,thin,dash pattern=on 0.025cm off 0.05cm]
        ({\cxi+\RRR2*cos(\aaa4-\dang2)},
         {\cyi+\RRR2*sin(\aaa4-\dang2)}) arc(\aaa4-\dang2:\aaa4-\dang2-\aaaa2:\RRR2);
    \draw[black,->,thin,dash pattern=on 0.025cm off 0.05cm] (\cxi,\cyi) -- +({\RRR1*cos(11)},{\RRR1*sin(11)}) node[anchor=west]{$\;\rho_\textit{min}$};
    \draw[black,->,thin,dash pattern=on 0.025cm off 0.05cm] (\cxi,\cyi) -- +({\RRR2*cos(4.5)},{\RRR2*sin(4.5)}) node[anchor=west]{$\;\rho_\textit{max}$};
    
    \draw[color4,line width=0.5mm]
       (\cx1,\cy1) to[out=160,in=0] (\cx2,\cy2) to[out=180,in=30] (\cx3,\cy3) to[out=210,in=115] (\cxi,\cyi);
    
	\coneback{\cx1}{\cy1}{\aaac1}{\Lcone}{\Acone}{\Dcone}{color=black!100,densely dotted}
    \coneback{\cx2}{\cy2}{\aaac2}{\Lcone}{\Acone}{\Dcone}{color=black!100,densely dotted}
    \coneback{\cx3}{\cy3}{\aaac3}{\Lcone}{\Acone}{\Dcone}{color=black!100,densely dotted}
	\draw[black,thin,->] (\cx2,\cy2) -- +({\pfact*\LL*sin(\aaac2)},{-\pfact*\LL*cos(\aaac2)});
	\draw[black] ({\cx2+\pfact*\LL*sin(\aaac2)},{\cy2-\pfact*\LL*cos(\aaac2)+0.15}) node[anchor=east] {$\losvec$};
    
    \rocket{\cx1}{\cy1}{\aa1}{\lrocket}{1.5*\lflame}{5}
    \rocket{\cx2}{\cy2}{\aa2}{\lrocket}{1.0*\lflame}{-5}
    \rocket{\cx3}{\cy3}{\aa3}{\lrocket}{1.5*\lflame}{-20}
    
    \draw[color0!100,line width=0.25mm] (\cx2,\cy2) -- (\cxi,\cyi);
    
    \cone{\cx1}{\cy1}{\aaac1}{\Lcone}{\Acone}{\Dcone}{color=color2!100,opacity=0.20,shading=axis,shading angle=90,left color=color2!100,right color=color2!0};
    \cone{\cx2}{\cy2}{\aaac2}{\Lcone}{\Acone}{\Dcone}{color=color1!100,opacity=0.20,shading=axis,shading angle=90,left color=color1!100,right color=color1!0};
    \cone{\cx3}{\cy3}{\aaac3}{\Lcone}{\Acone}{\Dcone}{color=color2!100,opacity=0.20,shading=axis,shading angle=90,left color=color2!100,right color=color2!0};
    
    \isoaxes{\cx2}{\cy2}{90+\aa2}{1.5}{12}{$\Fbody$}{$\bm{e}_1$}{$\bm{e}_2$}{$\bm{e}_3$}
    \isoaxes{\cxi}{\cyi}{90}{1.0}{15}{$\Finertial$}{}{}{}
    \fill[black] (\cx1,\cy1) circle (2pt);
    \fill[black] (\cx3,\cy3) circle (2pt);
    
    \draw[color=black!100,<-] ({\cx2+\fact*\LL*sin(\aaac2)},
                               {\cy2-\fact*\LL*cos(\aaac2)}) arc(-90+\aaac2:-90+\aaac2-\da:\fact*\LL); 
    \draw[color=black!100,<-] ({\cx2-\fact*\LL*cos(\aad)},
                               {\cy2-\fact*\LL*sin(\aad)}) arc(180+\aad:180+\aad+\da:\fact*\LL) node[anchor=west] {$\los$};
    
    \draw (\cxi,\cyi) to[out=-90,in=-180] (\cxi+1.0,\cyi-1.00) node[anchor=west] {Landing site};
    \draw (\cx2-2.8,\cy2+0.3) to[out=180,in=0] (\cx2-4.1,\cy2+1.3) node[anchor=east] {$\begin{aligned} \\ \\[0.5ex] \text{Trigger set where \LOS { cone}} \\ \text{constraint is activated} \\ \big\{\bm{r}:\rho_\textit{min} \leq \|\bm{r}\|_2 \leq \rho_\textit{max}\big\}\end{aligned}$};
    \draw (\cx3,\cy3) to[out=180,in=0] (\cx3-0.75,\cy3+0.50) node[anchor=east] {Loss of \LOS};
    \draw (\cx1,\cy1) to[out=0,in=180] (\cx1+0.5,\cy1-0.4) node[anchor=west] {PDI};
    \draw ({\cx1+\tfact*\LL*sin(\aaac1+0.5*\Acone)},
           {\cy1-\tfact*\LL*cos(\aaac1+0.5*\Acone)}) to[out=\aaac1+0.5*\Acone,in=180] (\cx1+0.5,\cy1-1.7) node[anchor=west] {$\begin{aligned} \\ &\text{Body-fixed}\\[-0.8ex] &\text{\LOS { cone}}\end{aligned}$};
\end{tikzpicture}
    \caption{Depiction of a slant-range-triggered line of sight (LoS) compound state-triggered constraint. The \LOS \ constraint is enforced only when the vehicle has a slant range between~$\rho_\textit{min}$ and~$\rho_\textit{max}$ to the landing site.}
    \label{fig:cstc_fov}
\end{figure}
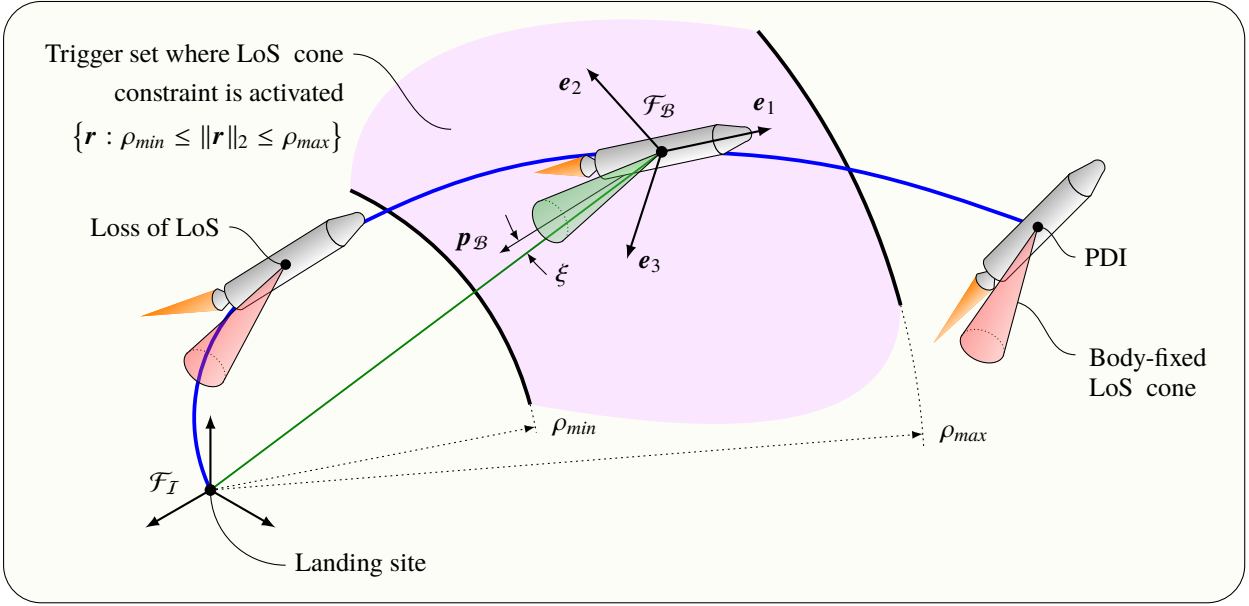

In this section we apply the compound state-triggered constraint given in~\eqref{eq:cstc_logic} and~\eqref{eq:cstc_cont} to a slant-range-triggered \LOS constraint. A scenario employing this constraint is illustrated in Fig.~\ref{fig:cstc_fov}, where the cone attached to the vehicle at each instance represents the field of view of an optical sensor. 
The trigger condition is given in inertial coordinates as
\begin{equation*}
    \bigwedge_{i=1}^{2}\big(\stcg_i(\rI)< 0\big) \definedas \big(\rho_\textit{min} < \|\rI\|_2\big)\land\big(\|\rI\|_2 < \rho_\textit{max}\big)
\end{equation*}
and is designed to trigger when the vehicle is at a range between~$\rho_\textit{min}$ and~$\rho_\textit{max}$ to the landing site. Expressing the trigger functions in terms of~$\dqq$ by using~\eqref{eq:norma_equiv_dq} gives
\begin{subequations} \label{eq:cstc_ex_trig}
    \begin{align}
        \stcg_1(\dqq) &\definedas \rho_\textit{min} - \|2 E_d \dqq\|_2, \\
        \stcg_2(\dqq) &\definedas \|2 E_d \dqq\|_2 - \rho_\textit{max}.
    \end{align}
\end{subequations}
The constraint condition limits the angle between $\losvec$ and the \LOS vector to the landing site, $-\rB$, to be less than~$\losmax\in(0\dg,90\dg)$ and is given by
\begin{equation}
    \rB\tran \losvec + \| \rB \|_2 \cos \losmax \leq 0.
    \label{eq:los_constraint_defn}
\end{equation}
For simplicity, we do not account for sensor offsets from the origin of the body frame. The left-hand side of~\eqref{eq:los_constraint_defn} can be expressed as a function of the dual quaternion $\dqq$ by using~\eqref{eq:quat_id_BB} and~\eqref{eq:norma_equiv_dq} to write
\begin{equation} \label{eq:cstc_ex_cons}
    \stcc(\dqq) \definedas \dqq\tran M_l \dqq + \|2E_d\dqq\|_2\cos\losmax, \quad M_l=\begin{bmatrix} 0_{4\times 4} &{\qskewstar{\losvec}}\tran \\ \qskewstar{\losvec} & 0_{4\times 4} \end{bmatrix}.
\end{equation}
References~\cite{Lee2015a,Lee2017} showed that the function~\eqref{eq:cstc_ex_cons} is a convex function of $\dqq\in\setudqR$ over the same domain as~\eqref{eq:dq_gs_max}. The compound state-triggered constraint is thus obtained by using~\eqref{eq:cstc_ex_trig} and~\eqref{eq:cstc_ex_cons} in the continuous formulation~\eqref{eq:cstc_cont}. 

\subsection{Non-convex Optimal Control Problem Statement}\label{sec2:ncvx_problem_statement}

We conclude this section by collecting the results of~\sref{sec2:baseline} and~\sref{sec2:stc} into a complete statement of the continuous-time non-convex optimal control problem to be solved. We focus on \textit{minimum fuel} problems, which are equivalently cast as maximizing the final mass. The problem statement is given in Problem~\ref{prob:prob1} where, for brevity, we have omitted the temporal arguments for all constraints except the boundary conditions.

\begin{problem}\label{prob:prob1}
Find the burn time, $\tf\in\realpp$, and the piecewise continuous thrust commands $\uB(t)$ over $t\in[\ti,\tf]$ that solve the optimal control problem:
\begin{equation*}
\begin{alignedat}{4}
\min_{t_f,\,\uB(\cdot)} &\quad &&-m(t_f) && \quad&&\text{convex} \\
\text{subject to:} &\quad &&\dot{m} = - \alpha \| \uB \|_2,\quad \dqqdot = \frac{1}{2} \dqq \otimes \dqw \qquad&&\eqref{eq:mass_dyn}\,\&\,\eqref{eq:dq_kinematics} \quad&&\text{non-convex} \\
&\quad &&\J \dqwdot = \Phi \uB + \dqgB - \dqw \oslash \J \dqw \qquad&&\eqref{eq:dq_dynamics} \quad&&\text{non-convex} \\
&\quad &&\gsf(\dqq) \leq 0, \quad \tiltf(\dqq) \leq 0, \quad \| \wB \|_{\infty} \leq \wmax, \qquad&&\eqref{eq:dq_gs_max},\,\eqref{eq:dq_tilt_max}\,\&\,\eqref{eq:w_max} \quad&&\text{convex} \\
&\quad &&\umin \leq \| \uB \|_2 \leq \umax \qquad&&\eqref{eq:thrust_bound} \quad&&\text{non-convex} \\
&\quad &&\| \uB \|_2 \leq \sec\gimbalmax \bm{z}_{\body}\tran \uB  \qquad&&\eqref{eq:max_gimbal} \quad&&\text{convex} \\
&\quad &&\|E_{xy}\duB\|_2 \leq \dgimmax \zB\tran\uB, \qquad&&\eqref{eq:max_gimbal_rate} \quad&&\text{convex} \\
&\quad &&-\dot{u}_{z,\max} \leq \zB\tran \duB \leq \dot{u}_{z,\max}, \qquad &&\eqref{eq:max_throttle_rate} \quad&&\text{convex}\\
&\quad &&\stch_{\land}(\dqq) \leq 0 \qquad&&\eqref{eq:cstc_ex_trig}\,\&\,\eqref{eq:cstc_ex_cons} \quad&&\text{non-convex} \\
&\quad &&m(\ti) = \mwet, \quad \dqq(\ti)=\dqqig(\q(\ti)), \quad \dqw(\ti) = \dqwig(\q(\ti)) \qquad&&\eqref{eq:bcs_init} \quad&&\text{non-convex} \\
&\quad &&\dqq(\tf) = \dqq_f, \quad \dqw(\tf)=\dqw_f , \quad m(\tf) \geq \mdry. \qquad&&\eqref{eq:bcs_final}\,\&\,\eqref{eq:max_fuel} \quad&&\text{convex}
\end{alignedat}
\end{equation*}
\end{problem}

\section{Solution Method}\label{sec:sol_method}

This section details the steps taken to solve Problem~\ref{prob:prob1}. First, a high-level overview of the algorithm is given in~\sref{sec2:successive_algo} and is depicted in Fig.~\ref{fig:scvx_loop}. The goal is to the transcribe Problem~\ref{prob:prob1} into a discrete-time parameter optimization problem that can be solved by using convex optimization solvers. Next, the three primary steps that are required to obtain a fully convex and numerically well-conditioned relaxation of Problem~\ref{prob:prob1} are introduced in~\sref{sec2:propagation_step},~\sref{sec2:parameter_update} and~\sref{sec2:solve_step}. A complete statement of the convexified problem is given in~\sref{sec2:solve_step} along with a discussion of the algorithm's stopping criterion.

\subsection{Successive Convexification Algorithm}\label{sec2:successive_algo}

\begin{figure}[b!]
    \centering
    \input{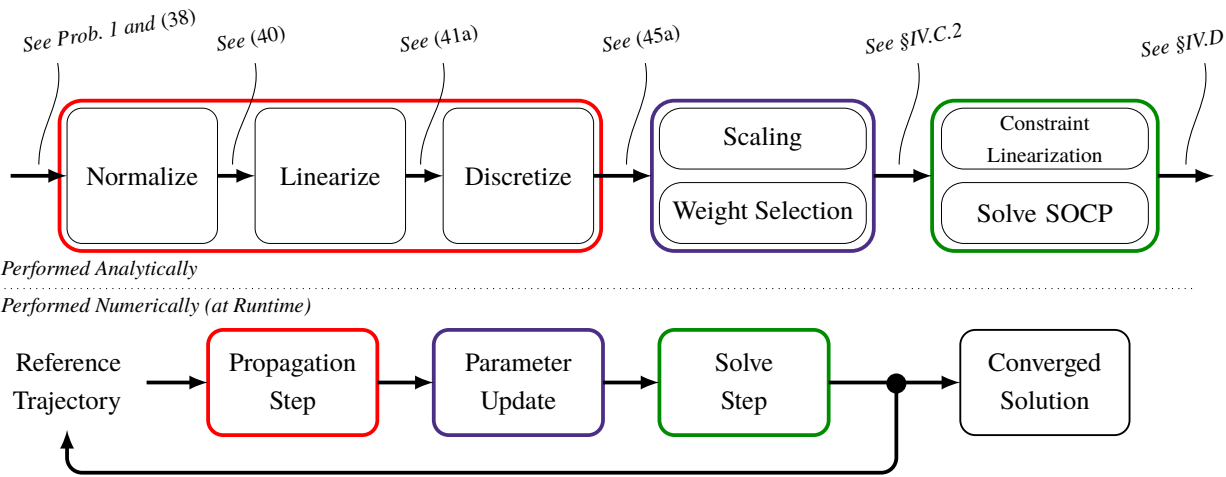}
    \caption{Overview of the successive convexification algorithm designed to solve Problem~\ref{prob:prob1}. The propagation step is detailed in~\sref{sec2:propagation_step}, the parameter update is detailed in~\sref{sec2:parameter_update} and the solve step is detailed in~\sref{sec2:solve_step}.}
    \label{fig:scvx_loop}
\end{figure}

Similar to~\cite{SzmukReynolds2018}, we seek a solution to Problem~\ref{prob:prob1} that approximates optimality, exactly satisfies the nonlinear equations of motion, and approximates the state and control constraints by enforcing them at a finite number of temporal nodes. The solution algorithm iteratively solves a sequence of subproblems, each of which represent a convex relaxation of Problem~\ref{prob:prob1}, and terminates once a convergence criterion is met. Each iteration consists of three primary steps: a \textit{propagation} step responsible for approximating the dynamic equations of motion, a \textit{parameter update} step responsible for maintaining proper numerical conditioning, and a \textit{solve} step responsible for solving each convex subproblem to full optimality. Together, these steps result in a burn time and corresponding state and control trajectory, or \textit{iterate}. The previous iterate is used as the basis for computing all approximations used during the current iteration. In this work, the parameter update step extends this relationship by directly using information from the current iteration's propagation step to better inform the \textit{current iteration's} solve step. This modification has improved convergence in numerical trials compared to those reported in~\cite{SzmukReynolds2018,Szmuk2019,Reynolds2019,Malyuta2019}. The bottom half of Fig.~\ref{fig:scvx_loop} provides a block diagram representation of the algorithm's operation, while the top half shows the same steps broken out into the analytical operations discussed in~\sref{sec2:propagation_step}-\sref{sec2:solve_step}.

\subsection{Propagation Step}\label{sec2:propagation_step}

Equality constraints for a general convex optimization problem must be affine functions of the solution variables. As such, the goal of this section is to convert the continuous-time nonlinear equations of motion~\eqref{eq:mass_dyn},~\eqref{eq:dq_kinematics} and~\eqref{eq:dq_dynamics} into discrete-time affine functions of the state, control, and burn time. To outline the method, we use a general autonomous nonlinear system of the form
\begin{equation}
  \label{eq:general_nonlinear_dynamics}
  \dot \xx(t) = f(\xx(t),\uu(t)),\quad\forall t\in [\ti,\tf],
\end{equation}
where we assume that $\xx(t)\in\real^{n_x}$ and $\uu(t)\in\real^{n_u}$. The approach closely follows \cite{SzmukReynolds2018,Malyuta2019} but is summarized here for the sake of completeness. To begin, the dynamics are temporally normalized to convert~\eqref{eq:general_nonlinear_dynamics} into an equivalent fixed-final-time expression. We define the \textit{scaled time} $\tau\in [0,1]$ such
that
\begin{equation}
  \label{eq:real_scaled_time_relation}
  t = \ss\tau,
\end{equation}
where $s>0$ is a temporal dilation factor. An application of the chain rule yields
\begin{equation}
  \label{eq:normalized_nonlinear_dynamics}
  \xx'(\tau) \definedas \frac{d}{d\tau}\xx(\tau)
  = \frac{dt}{d\tau}\frac{d}{dt}\xx(\tau)
  = \ss f(\xx(\tau),\uu(\tau)) \definedas F(\ss,\xx(\tau),\uu(\tau)).
\end{equation}

Next,~\eqref{eq:normalized_nonlinear_dynamics} must be linearized in order to be an affine function of the solution variables. We approximate the dynamics with a first-order Taylor series expansion about a time-varying reference trajectory
$\zzb(\tau)~\definedas~[\ssb ~\xxb(\tau)^T ~\uub(\tau)^T ]^T$ for all
$\tau\in [0,1]$. Section~\ref{sec3:linearization_reference_trajectory}
details how this reference trajectory is obtained from a previous iterate. A Taylor series expansion of~\eqref{eq:normalized_nonlinear_dynamics} leads to a linear time-varying (LTV) system:
\begin{subequations}
  \begin{align}
    \label{eq:linearized_continuous_dynamics}
    \xx'(\tau)
    &\approx A(\tau)\xx(\tau)+B(\tau)\uu(\tau)+S(\tau)\ss+\ww(\tau),
      \quad\forall\tau\in [0,1], \\
    A(\tau)
    &\definedas \left.\frac{\partial F}{\partial \xx}\right|_{\zzb (\tau)}, \\
    B(\tau)
    &\definedas \left.\frac{\partial F}{\partial \uu}\right|_{\zzb (\tau)}, \\
    S(\tau)
    &\definedas \left.\frac{\partial F}{\partial \ss}\right|_{\zzb (\tau)}, \\
    \ww (\tau)
    &\definedas -A(\tau)\xxb(\tau) - B(\tau) \uub(\tau).
  \end{align}
\end{subequations}

The LTV system~\eqref{eq:linearized_continuous_dynamics} is discretized by representing the infinite dimensional control signal $\uu(\tau)$ with a finite number of parameters. In this sense, our method constitutes a direct method for solving optimal control problems~\cite{Betts1998}.  We divide the scaled time into $N-1$ subintervals and approximate the control signal as a continuous piecewise-affine function by using the so-called affine interpolation
\begin{subequations}
  \begin{align}
    \label{eq:foh_input}
    \uu(\tau) &= \lambda_\ind^-(\tau)\uu_{\ind}+\lambda_{\ind}^+(\tau)\uu_{\ind+1},
    \quad\forall\tau\in[\tau_{\ind},\tau_{\ind+1}], \\
    \lambda_{\ind}^-(\tau)&\definedas\frac{\tau_{\ind+1}-\tau}{\tau_{\ind+1}-\tau_{\ind}},\quad
    \lambda_{\ind}^+(\tau)\definedas\frac{\tau-\tau_{\ind}}{\tau_{\ind+1}-\tau_{\ind}},
  \end{align}
\end{subequations}
for each $\ind\in\setKb$, where $\tau_{\ind}=\nicefrac{\ind-1}{N-1}$. The $\uu_{\ind}$ are now constant parameters for each $\ind\in\setK$ that can be used to reconstruct the continuous-time signal $\uu(\tau)$ for any $\tau\in[0,1]$. Using~\eqref{eq:foh_input} in~\eqref{eq:linearized_continuous_dynamics}, the LTV dynamics over each time interval become
\begin{equation}
  \label{eq:continuous_foh_dynamics}
  \xx'(\tau) = A(\tau)\xx(\tau)+\lambda_{\ind}^-(\tau)B(\tau)\uu_{\ind}+\lambda_{\ind}^+(\tau)B(\tau)\uu_{\ind+1}+
  \SS(\tau)\ss+\ww(\tau),\quad\tau\in [\tau_{\ind},\tau_{\ind+1}].
\end{equation}

The zero-input state transition matrix
$\Phi(\cdot,\tau_\ind):[\tau_{\ind},\tau_{\ind+1}]\to\real^{n_x\times n_x}$ associated with
\eqref{eq:continuous_foh_dynamics} is given by
\begin{equation}
  \label{eq:stm}
  \Phi(\tau,\tau_{\ind}) = I + \int_{\tau_{\ind}}^\tau A(\dummy)\Phi(\dummy,\tau_{\ind})d\dummy.
\end{equation}

Using the inverse and transitive properties of the state transition matrix~\cite{Antsaklis2007}, a discrete-time LTV equation can be obtained for each $\ind\in\setKb$:
\begin{subequations} 
\begin{align}
\xx_{\ind+1} &= A_{\ind} \xx_{\ind} + B_{\ind}^-\uu_{\ind} + B_{\ind}^+\uu_{\ind+1} + S_{\ind} \ss + \ww_{\ind}, \label{eq:dltv} \\
A_{\ind} &\definedas \Phi(\tau_{\ind+1},\tau_{\ind}), \label{eq:dltv_A} \\
B_{\ind}^- &\definedas A_{\ind} \int_{\tau_{\ind}}^{\tau_{\ind+1}} \Phi^{-1}(\tau,\tau_{\ind})\lambda_{\ind}^-(\tau)B(\tau)d\tau, \label{eq:dltv_Bm} \\
B_{\ind}^+ &\definedas A_{\ind}\int_{\tau_{\ind}}^{\tau_{\ind+1}} \Phi^{-1}(\tau,\tau_{\ind})\lambda_{\ind}^+(\tau)B(\tau)d\tau, \label{eq:dltv_Bp} \\
\SS_{\ind} &\definedas A_{\ind} \int_{\tau_{\ind}}^{\tau_{\ind+1}} \Phi^{-1}(\tau,\tau_{\ind})\SS(\tau)d\tau,  \label{eq:dltv_S}\\
\ww_{\ind} &\definedas A_{\ind} \int_{\tau_{\ind}}^{\tau_{\ind+1}} \Phi^{-1}(\tau,\tau_{\ind}) \ww (\tau)d\tau. \label{eq:dltv_w}
  \end{align}
\end{subequations}

In implementation, \eqref{eq:stm} and \eqref{eq:dltv_Bm}-\eqref{eq:dltv_w} are integrated simultaneously via the classical Runge-Kutta RK4 method~\cite{Butcher2008}. 
The final results are obtained after the numerical integration by right multiplying the final value of~\eqref{eq:stm} by the final value of each integral in~\eqref{eq:dltv_Bm}-\eqref{eq:dltv_w}.

\subsubsection{Reference Trajectory Computation}\label{sec3:linearization_reference_trajectory}

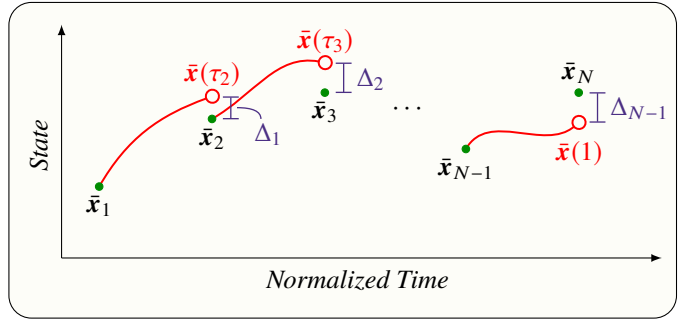
\begin{figure}[t!]
  \centering




\begin{tikzpicture}[shorten >=-3pt]
    \newcommand{\circweight}{0.05}
    \newcommand{\lenw}{0.75pt}
    
    \tikzmath{\dt=1.5; \dx=0.25; \dy=0.025;
              \ddx=0.4; \ddy=0.075;
              \pppx1=0.5;       \pppy1=0.95;
              \pppx2=\pppx1+\dt; \pppy2=\pppy1+0.90;
              \pppx3=\pppx2+\dt; \pppy3=\pppy2+0.35;
              \pppx4=\pppx3+0.75*\dt; \pppy4=\pppy3-0.25;
              \pppx5=\pppx4+0.5*\dt; \pppy5=\pppy4-0.50;
              \pppx6=\pppx5+\dt; \pppy6=\pppy5+0.75;
              \lx=\pppx6+1; \ly=3.0; 
              \ppx1=0.00; \ppy1=2.2;
              \ppx2=\lx-0.2; \ppy2=0.5;
    }
    
    \pane{0.5*\lx-0.2}{0.5*\ly-0.2}{0}{\lx+1.0}{\ly+1.2}{0.3}{color=black,fill=beige!20};
    
    \draw[black,->] (0,0) --+(\lx,0);
    \draw (0.5*\lx,-0.3) node[anchor=center] {\it Normalized Time};
    
    \draw[black,->] (0,0) --+(0,\ly);
    \draw (-0.3,0.5*\ly) node[anchor=center,rotate=90] {\it State};
    
    \draw[propcolor,line width=\lenw,cap=round,-o] (\pppx1,\pppy1) to[out=60,in=200] (\pppx2,\pppy2+0.3) node[anchor=south] {$\xxb(\tau_2)$};
    \draw[propcolor,line width=\lenw,cap=round,-o] (\pppx2,\pppy2) to[out=30,in=170] (\pppx3,\pppy3+0.4) node[anchor=south] {$\xxb(\tau_3)$};
    \draw[propcolor,line width=\lenw,cap=round,-o] (\pppx5,\pppy5) to[out=60,in=220] (\pppx6,\pppy6-0.4) node[shift={(0,-0.4)}] {$\xxb(1)$};
    
    \filldraw[solvecolor] (\pppx1,\pppy1) circle (\circweight);
    \draw (\pppx1,\pppy1) node[anchor=north] {$\xxb_1$};
    
    \filldraw[solvecolor] (\pppx2,\pppy2) circle (\circweight);
    \draw (\pppx2,\pppy2) node[anchor=north] {$\xxb_2$};
    
    \filldraw[solvecolor] (\pppx3,\pppy3) circle (\circweight);
    \draw (\pppx3,\pppy3) node[anchor=north] {$\xxb_3$};
    
    \draw (\pppx4,\pppy4) node {$\cdots$};

    \filldraw[solvecolor] (\pppx5,\pppy5) circle (\circweight);
    \draw (\pppx5,\pppy5) node[anchor=north] {$\xxb_{N-1}$};

    \filldraw[solvecolor] (\pppx6,\pppy6) circle (\circweight);
    \draw (\pppx6,\pppy6) node[anchor=south] {$\xxb_N$};
       
    \draw[scalecolor,line width=0.5pt,arrows=|-|] (\pppx2+\dx,\pppy2+0.3) -- (\pppx2+\dx,\pppy2+0.1);
    \node[scalecolor] at (\pppx2+0.75,\pppy2-0.2) {$\res_1$};
    \draw[scalecolor,line width=0.5pt] (\pppx2+0.65,\pppy2) to[out=90,in=0] (\pppx2+1.75*\dx,\pppy2+0.15);
    \draw[scalecolor,line width=0.5pt,arrows=|-|] (\pppx3+\dx,\pppy3+0.4) -- (\pppx3+\dx,\pppy3+0.1) node[shift={(0.9*\ddx,\ddy)}] {$\res_{2}$};
    \draw[scalecolor,line width=0.5pt,arrows=|-|] (\pppx6+\dx,\pppy6-0.4) -- (\pppx6+\dx,\pppy6-0.1) node[shift={(1.4*\ddx,-1.5*\ddy)}] {$\res_{N-1}$};

\end{tikzpicture}

  \caption{Illustration of the integration procedure in~\eqref{eq:reference_state_integration} used to obtain the reference state trajectory. The colors match those of Fig.~\ref{fig:scvx_loop}; the green dots are obtained from the previous solve step, and the continuous red curves serve as the reference trajectory for the current propagation step. The purple defects (see~\eqref{eq:defects}) represent the disagreement between the numerically integrated red curves and predicted green dots, and are used in the parameter update step.}
  \label{fig:integration_resetting}
\end{figure}

To obtain the discrete-time LTV system~\eqref{eq:dltv}, a reference trajectory $\zzb(\tau)$ was assumed to be available for any time point $\tau\in[\tau_{\ind},\tau_{\ind+1}]$ and for each $\ind\in\setKb$. For the first iteration, the reference trajectory is obtained by using the straight-line initialization method presented in~\cite{SzmukReynolds2018}. Subsequent reference trajectories are obtained by using the previously computed iterate, which provides $\ssb$, $\{\uub_{\ind}\}_{\ind=1}^N$ and $\{\xxb_{\ind}\}_{\ind=1}^N$. Using these discrete vectors, $\uub(\tau)$ is obtained for any $\tau\in[\tau_{\ind},\tau_{\ind+1}]$ by using $\uub_{\ind}$ and $\uub_{\ind+1}$ in accordance with~\eqref{eq:foh_input}. As illustrated in Fig.~\ref{fig:integration_resetting}, the reference state trajectory on this same time interval is obtained by integrating the continuous control signal through the nonlinear dynamics according to
\begin{equation}
  \label{eq:reference_state_integration}
  \xxb(\tau) = \xxb_{\ind}+\int_{\tau_{\ind}}^\tau  F(\ssb,\,\xxb(\dummy),\uub(\dummy))d\dummy,
  \quad\tau\in[\tau_{\ind},\tau_{\ind+1}].
\end{equation}

We note that \eqref{eq:reference_state_integration} \textit{resets} the state
trajectory to evolve from $\xxb_{\ind}$ over each subinterval, as shown in Fig.~\ref{fig:integration_resetting}. This is analogous to a multiple shooting technique, and -- similar to how that method exhibits improved convergence when compared to single shooting techniques -- we have observed that~\eqref{eq:reference_state_integration} improves convergence by keeping
$\xxb(\cdot)$ closer to the individual points $\{\xxb_{\ind}\}_{\ind=1}^N$~\cite{SzmukReynolds2018,Malyuta2019}. For later use, we define the \textit{defect} associated with the $\ind$th interval to be
\begin{equation}
    \resk \definedas \| \xxb(\tau_{\ind+1}) - \xxb_{\ind+1} \|, \quad \ind\in\setKb.
    \label{eq:defects}
\end{equation}
These defects are shown in Fig.~\ref{fig:integration_resetting} and compare the final result from~\eqref{eq:reference_state_integration} with the predicted values from the previous iterate at the same time node. 

\subsection{Parameter Update Step}\label{sec2:parameter_update}

Linearizing and discretizating the dynamics in~\sref{sec2:propagation_step} results in an a discrete-time affine representation of the nonlinear dynamics that can be used within a convex optimization problem. 
This section will first detail the use of \textit{trust regions} and \textit{virtual control} in the context of the successive convexification technique used in this work to guide the iterative process. Here, we present a new heuristic method for determining the trust region sizes that leverages information available from the propagation step. Next, we discuss a scaling method that helps to ensure a numerically well-conditioned optimization problem in the subsequent solve step. 

\subsubsection{Adaptive Trust Region Weighting and Virtual Control}\label{sec3:vc_and_tr}

Iterative algorithms based on the linearization of nonlinear or non-convex terms are generally subject to a trust region constraint to ensure that each iterate is chosen from a region where the linearization is valid~\cite{NocedalWright}. In the context of aerospace problems, both hard trust regions~\cite{Mao2016a,Mao2017,Mao2018,LiuSurvey2017} and soft trust regions~\cite{Szmuk2016,Szmuk2017,Szmuk2018,SzmukReynolds2018} have been applied successfully. While using hard trust regions may be more prudent \textit{in general} to use~\cite{NocedalWright}, we have found via extensive numerical experimentation that allowing the optimization process to select the trust-region radius for powered descent problems works well. We therefore augment our problem with the quadratic constraint
\begin{equation}
    \left\| \xx_{\ind} - \xxb_{\ind} \right\|_2^2 + \left\| \uu_{\ind} - \uub_{\ind} \right\|_2^2 \leq \eek, \quad \ind\in\setK,
    \label{eq:trust_region}
\end{equation}
where $\ee\in\real^N$ is a vector of trust region radii that are added as solution variables. To encourage shrinking trust region radii we augment the cost function with
\begin{equation}
    \Jtr = \wtr ^T \ee,
    \label{eq:trust_region_cost}
\end{equation}
where $\wtr\in\realpp^N$ is a vector of trust region weights. Note that the term in~\eqref{eq:trust_region_cost} is effectively a weighted 1-norm penalty on the trust region radii due to~\eqref{eq:trust_region}.  We use a new method to select the weights $\wtr$ such that they reflect the accuracy of the discretization routine in~\sref{sec2:propagation_step}. In particular, we use the defects computed in~\eqref{eq:defects} to compute the weights as
\begin{equation}
    \wtrk = \begin{cases}
    \frac{1}{\|\resk\|}, &\quad \|\resk\| \geq \resmin \\
    \frac{1}{\resmin}, &\quad \| \resk \| < \resmin
    \end{cases}
    \label{eq:tr_weight_update}
\end{equation}
for each $\ind\in\setKb$, where $\resmin$ provides an upper bound on the weighting terms that helps to maintain proper scaling and avoid numerical issues when the residual terms become small. The result is that if the propagation step indicates a small defect at the $\ind$th time node, then the cost of deviating from the previous iterate's values at this node is increased in the subsequent solve step. This guides the solution process towards feasibility with respect to the nonlinear dynamics, while still permitting variations for the sake of constraint satisfaction and optimality. 

The use of trust regions can create another issue when used in conjunction with linearization-based iterative methods; namely that of \textit{artificial infeasibility}~\cite{NocedalWright,Mao2016a,Mao2017,Mao2018}. This can arise when satisfaction of a linearized constraint outweighs the cost of satisfying the trust-region constraint -- i.e., the constraints are inconsistent and cannot be simultaneously satisfied without taking too large a step. To alleviate these issues we use the idea of \textit{virtual control} that acts as a synthetic and unconstrained input. In our implementation we limit the use of virtual control to the equality constraint that arises from the equations of motion~\eqref{eq:dltv}, which is augmented with the virtual control term $\vv_{\ind}\in\real^{n_x}$ according to
\begin{equation}
    \xx_{\ind+1} = A_{\ind} \xx_{\ind} + B_{\ind}^- \uu_{\ind} + B_{\ind}^ + \uu_{\ind+1} + \SS_{\ind} \ss + \ww_{\ind} + \vv_{\ind}, \quad \ind\in\setKb.
    \label{eq:dltv_with_vc}
\end{equation}
To penalize the use of virtual control such that it is only used when necessary for constraint satisfaction, we augment the cost with
\begin{equation}
    \Jvc = \wvc \sum_{\ind=1}^{N-1} \| \vv_\ind \|_1,
\end{equation}
where $\wvc\in\realpp$ is a large weighting term. Finally, it is important to note that when both $\Jvc$ and $\Jtr$ are zero, the augmented cost function is equivalent to that of Problem~\ref{prob:prob1}.

\subsubsection{Scaling}\label{sec3:scaling}

In theory, the performance and solution of an optimization problem are invariant to the magnitudes of the components of $\xx_{\ind}$ and $\uu_{\ind}$. In practice, however, numerical issues associated with sensitivity, machine precision and the choice of termination criteria can arise in a numerical solver when these magnitudes are highly different~\cite{NocedalWright,Gill1981}. A canonical example would be to optimize at once over thrust, which can have values on the order of $10^4~\unit{N}$, and unit quaternions, whose components are between $-1$ and $1$. We use the standard remedy of applying the affine variable transformations~\cite{NocedalWright,BettsBook}
\begin{subequations}\label{eq:scale}
  \begin{align}
    \xx_{\ind} &=\Px \xxs_{\ind}+\qx, \label{eq:scale_x} \\
    \uu_{\ind} &= \Pu \uus_{\ind}+\qu, \label{eq:scale_u} \\
    s &= \Pt \sss. \label{eq:scale_s}
  \end{align}
\end{subequations}
The scaling terms are chosen such that the components of the scaled state~$\xxs_{\ind}$, the scaled input~$\uus_{\ind}$ and the
scaled dilation factor $\sss$ have a maximum magnitude of roughly unity across all iterations. In particular, we choose:
\begin{subequations}
  \begin{align}
    \Px &= \diag{
          \mi,\,\ones{4},\,\tfrac{1}{2}\|\ri{\inertial}\|_2\ones{4},\,
          \wmax \ones{4},\,
          \|\vi{\inertial}\|_2\ones{4}
          }, \\
    \Pu &= \diag{\umax\ones{3}},
  \end{align}
\end{subequations}
where $\ones{p}\in\real^p$ is a vector of ones, and $\qx=0$ and $\qu=0$. For $\Pt$, we choose the optimal value of $\ssb$ from the algorithm's previous iterate. We note that other choices are possible for the scaling terms, such as to balance the magnitudes of the dual variables~\cite{Ross2018}, to minimize the condition number of the cost function Hessian at the solution, or to improve the behavior of the cost function's first and second derivatives with respect to machine precision~\cite{Gill1981}. In any case, \eqref{eq:scale_x}-\eqref{eq:scale_s} are used to replace $\xx_{\ind}$, $\uu_{\ind}$ and $\ss$ everywhere in the
optimization problem.

\begin{remark}
To ensure proper scaling of the auxiliary variable $\ee$, the trust region~\eqref{eq:trust_region} should be implemented using $\| \xxs_{\ind} - \Px^{-1}(\xxb_{\ind}-\qx)\|_2^2 + \| \uus_{\ind} - \Pu^{-1}(\uub_{\ind}-\qu)\|_2^2 \leq \eeks$. The variable $\ees\ne\ee$ represents the trust-region radii for the scaled problem. The virtual control terms will be naturally scaled and do not need to be treated further.
\end{remark}

\subsection{Solve Step}\label{sec2:solve_step}

This section summarizes the convex subproblem that is solved to full optimality during each solve step. Each convex subproblem is solved using a customizable interior point method algorithm designed for convex second-order cone programs~\cite{Dueri2014,Dueri2016}. Specific details of the solver are not addressed in this paper. 
A discussion of the overall algorithm's stopping criteria is provided in~\sref{sec3:convergence}.
%
%
While~\sref{sec2:propagation_step} outlined how the nonlinear dynamical equations are handled, two sources of nonconvexity remain in our problem. These come from the boundary conditions in~\eqref{eq:bcs_init} and the state-triggered constraint introduced in~\eqref{eq:cstc_ex_trig}-\eqref{eq:cstc_ex_cons}. Consider the following generalized non-convex constraint that is applied at some $\ind\in\setK$
\begin{equation}
  \label{eq:general_ncvx_constraint}
  h(\bm{z}_{\ind}) \le 0,
\end{equation}
where $h(\cdot):\real^{n_z}\to\real$ is assumed to be at least once differentiable almost everywhere. We approximate the constraint with its first-order Taylor series expansion about the reference trajectory from~\sref{sec3:linearization_reference_trajectory} as
\begin{equation}
  h(\zzb_{\ind})+\left.\frac{\mathrm{d}h}{\mathrm{d}\zz}\right|_{\zzb_{\ind}}
  \left( \bm{z}_{\ind}-\bar{\bm{z}}_{\ind} \right) \leq 0.
  \label{eq:general_linearized_constraint}
\end{equation}
The equality constraints~\eqref{eq:bcs_init} are treated similarly. Since $\dqqig$ is a linear function of the initial attitude, only $\dqwig$ in~\eqref{eq:bcs_init} is enforced using the analogue of~\eqref{eq:general_linearized_constraint}.

\subsubsection{Second-Order Cone Relaxation of Problem~\ref{prob:prob1}}

We are now prepared to state the convex parameter optimization problem that is a relaxation of Problem~\ref{prob:prob1}. Recall that $\xx_\ind = \big[ m_\ind \ \dqq_\ind \ \dqw_\ind \big]^T$, $\uu_\ind = \uk{\body}$ and $\ss = \tf$, and wherever they appear are equivalent to $\xx_{\ind} = \Px \xxs_{\ind} + \qx$, $\uu_{\ind} = \Pu \uus_{\ind} + \qu$ and $\ss = \Pt \sss$. We define $\bm{c} \definedas [-1 \ \zeros{16}{1} ]^T$ and $\dtau \definedas \tau_{\ind+1}-\tau_{\ind}$ for any $\ind\in\setKb$. The convex relaxation of Problem~\ref{prob:prob1} is given by Problem~\ref{prob:prob2}.

\begin{problem}\label{prob:prob2}
Let $\XX = \{ \xx_\ind \}_{\ind=1}^{N}$, $\UU=\{\uu_\ind\}_{\ind=1}^{N}$ and $\VV=\{\vv_\ind\}_{\ind=1}^{N-1}$ represent the composite state, control and virtual control vectors, respectively. Their scaled analogues are $\XXs,\,\UUs$ and $\VVs=\VV$. Find $\XXs\in\real^{Nn_x},\,\UUs\in\real^{N n_u},\,\sss\in\realpp,\,\VVs\in\real^{(N-1) n_x}$ and $\ees\in\real^N$ that solve the following parameter optimization problem: 
\begin{equation*}
\begin{alignedat}{3}
\min_{\XXs,\UUs,\sss,\VVs,\ees} &\quad && \bm{c}^T \xxs_N + \Jtr + \Jvc & \\
\text{subject to:} &\quad &&\xx_{\ind+1} = A_{\ind} \xx_{\ind} + B_{\ind}^- \uu_{\ind} + B_{\ind}^+ \uu_{\ind+1} + \SS_{\ind} \ss + \ww_{\ind} + \vv_{\ind}, \qquad&&\eqref{eq:dltv_with_vc} \\
&\quad &&\| \xxs_{\ind} - \Px^{-1}(\xxb_{\ind}-\qx) \|_2 + \| \uus_{\ind} - \Pu^{-1}(\uub_{\ind}-\qu) \|_2 \leq \eeks, \qquad&&\eqref{eq:trust_region} \\
&\quad &&\gsf(\dqq_{\ind}) \leq 0, \quad \tiltf(\dqq_{\ind}) \leq 0, \quad \| H_w \dqw_{\ind} \|_{\infty} \leq \wmax, \qquad&&\eqref{eq:dq_gs_max},\,\eqref{eq:dq_tilt_max}\,\&\,\eqref{eq:w_max} \\
&\quad &&\| \uu_{\ind} \|_2 \leq \umax, \quad \umin - \frac{\uub_{\ind}^T}{\| \uub_{\ind}\|_2} \uu_{\ind} \leq 0, \qquad&&\eqref{eq:thrust_bound} \\
&\quad &&\| \uu_{\ind} \|_2 \leq \sec \gimbalmax \, \zB^T \uu_{\ind},  \qquad&&\eqref{eq:max_gimbal} \\
&\quad &&\|E_{xy}\left(\nicefrac{\uu_{\ind+1} - \uu_{\ind}}{\dtau}\right)\|_2 \leq \dgimmax \zB^T \uu_{\ind}, \qquad&&\eqref{eq:max_gimbal_rate} \\
&\quad &&-\dot{u}_{z,\max}\dtau \leq \zB ^T \left(\uu_{\ind+1}-\uu_{\ind}\right) \leq \dot{u}_{z,\max}\dtau, \qquad&&\eqref{eq:max_throttle_rate}\\
&\quad &&\stch_{\land}(\dqb_{\ind}) + \left.\frac{\partial \stch_{\land}}{\partial \dqq}\right|_{\dqb_{\ind}}^T (\dqq_{\ind} - \dqb_{\ind}) \leq 0, \qquad&&\eqref{eq:cstc_ex_trig}\,\&\,\eqref{eq:cstc_ex_cons} \\
&\quad &&m_1 = \mi, \quad \dqq_1=\dqqig(\q_1), \quad \dqw_1 = \dqwig(\qb_1) + \left.\frac{\partial \dqwig}{\partial \q}\right|_{\qb_1} (\q_1 - \qb_1), \qquad&&\eqref{eq:bcs_init} \\
&\quad &&\dqq_N = \dqq_f, \quad \dqw_N =\dqw_f, \quad m_N \geq \mdry. \qquad&&\eqref{eq:bcs_final}\,\&\,\eqref{eq:max_fuel}
\end{alignedat}
\end{equation*}
\end{problem}

\subsubsection{Convergence and Stopping Criteria}\label{sec3:convergence}

An important consideration for any iterative algorithm is the definition of ``convergence'' that is used to decide when to stop the iterations. Typically this is achieved by comparing some notion of difference between the solutions of two subsequent iterates. Our approach is to compare the maximum deviation in the state solution between two iterations, and to terminate when this difference is less than a prescribed tolerance. We do not consider changes in the control signal since large deviations in thrust commands may cause negligible changes to the trajectory and can be an overly conservative indication of convergence. Moreover, changes in the state are best determined by using the \textit{scaled} values discussed in~\sref{sec3:scaling}. Observing changes in the scaled quantities ensures that the notion of ``small'' is uniform across all state variables. This ensures -- for example -- that we do not treat a $1\,\unit{m}$ deviation in position the same as a $1\,\unit{rad/s}$ change in angular velocity. Given some tolerance $\dxtol\in\realpp$, the iterations are terminated whenever
\begin{equation}
    \max_{\ind\in\setK} \| P_x^{-1} (\xx_{\ind} - \xxb_{\ind}) \|_{\infty} < \dxtol.
    \label{eq:convergence}
\end{equation}

\begin{remark}
The virtual control and trust region sizes can also be an effective measure of convergence. Near zero values for both $\| \VV \|_1$ and $\| \ee \|_1$ can indicate that it is appropriate to terminate at the current iterate. However, we have observed that different numerical solvers can report different values for these two quantities, even for the same problem. To avoid solver-specific termination criteria, the state-based condition in~\eqref{eq:convergence} is used and has been observed to be uniformly applicable across numerical solvers.
\end{remark}

\section{Numerical Case Studies}\label{sec:numerical_demo}

This section provides two numerical case studies that highlight the capabilities of the solution method presented in~\sref{sec:sol_method}. First, the slant-range-triggered \LOS constraint is studied in~\sref{sec2:stc_case_study}. Second, a Monte Carlo analysis is performed in~\sref{sec2:monte_carlo} to study the algorithm's properties. For each case, we use the same set of nominal vehicle parameters outlined in Table~\ref{tab:params_ex}. The lunar gravitational acceleration is assumed to be $\gI = -1.62\zI$ and the inertia matrix is approximated for an Apollo-class vehicle with the initial mass listed in Table~\ref{tab:params_ex}. Recall that the baseline problem includes all constraints from \sref{sec2:baseline} and is equivalent to Problem~\ref{prob:prob1} \textit{without} the slant-range-triggered \LOS constraint.

\newcommand{\tabspace}{\hspace{-0.35cm}}
\newcommand{\tabsep}{\null\hspace{0.5cm}\null}
\begin{table}[b!]
  \centering
  \caption{Problem parameters for state-triggered constraint case study and baseline problem for the Monte Carlo case study.}
  \begin{tabular}{lrllclrll}
  	\hhline{====~====}
    \textbf{Parameter}      &   & \tabspace\textbf{Value}                   & \textbf{Units} &\tabsep&
    \textbf{Parameter}      &   & \tabspace\textbf{Value}                   & \textbf{Units} \\
    \hhline{----~----}
	$\mi$                   &   & \tabspace$3250$                           & $\unit{kg}$ &\tabsep&
	$\mdry$                 &   & \tabspace$2100$                           & $\unit{kg}$ \\
	$\ri{\inertial}$        &   & \tabspace$[250\;0\;433]$                  & $\unit{m}$ &\tabsep& 
	$\vi{\inertial}$        &   & \tabspace$[-30\;0\;-15]$                  & $\unit{m/s}$ \\
	$\rf{\inertial}$        &   & \tabspace$[0\;0\;30]$                     & $\unit{m}$ &\tabsep&
	$\vf{\inertial}$        &   & \tabspace$[0\;0\;-1]$                     & $\unit{m/s}$ \\
    $\ru$                   &$-$& \tabspace$0.25\cdot\zB$                   & $\unit{m}$ &\tabsep&
	$\wi,\,\wf$             &   & \tabspace$[0\;0\;0]$                      & $\unit{\dg/s}$ \\
    $I_{\text{sp}}$         &   & \tabspace$225.0$                          & $\unit{s}$ &\tabsep&
	$\qf$                   &   & \tabspace$[0\;0\;0\;1]$                   & - \\
    $\tiltmax$              &   & \tabspace$80.0$                           & $\unit{\dg}$ &\tabsep&
	$\wmax$                 &   & \tabspace$28.6$                           & $\unit{\dg/s}$ \\
	$\gsmax$                &   & \tabspace$80.0$                           & $\unit{\dg}$ &\tabsep&
	$\gimbalmax$            &   & \tabspace$20.0$                           & $\unit{\dg}$ \\
    $\umin$                 &   & \tabspace$6000$                           & $\unit{N}$ &\tabsep&
    $\umax$                 &   & \tabspace$22,500$                         & $\unit{N}$ \\
    $\dgimmax$               &   & \tabspace$n/a$                           & $\unit{rad/s}$ &\tabsep&
    $\dot{u}_{z,\max}$                 &   & \tabspace$n/a$                         & $\unit{N/s}$ \\
	\hhline{====~====}
    \label{tab:params_ex}
  \end{tabular}
\end{table}
\vspace{-12pt}
\subsection{State-Triggered Constraint Case Study}\label{sec2:stc_case_study}

This case study examines how the state-triggered constraint introduced in~\sref{sec3:stc_los} can be used to solve for powered descent trajectories when attitude pointing requirements must be met only during certain portions of the descent. In addition to the parameters in Table~\ref{tab:params_ex}, we assume in this section that an optical sensor has a boresight vector $\losvec = [ 0.91\; 0\; -0.42 ]^T$ in the body frame. Given current sensor technology under study for future missions~\cite{DwyerCianciolo2019}, we assume a field of view angle of $\losmax = 20\dg$. We impose the \LOS constraint whenever the vehicle's slant range lies between $\rho_\textit{min} = 200\unit{m}$ and $\rho_\textit{max} = 450\unit{m}$. 

\begin{figure}[tb]
    \centering
    \subfloat[]{\includegraphics[width=0.48\textwidth]{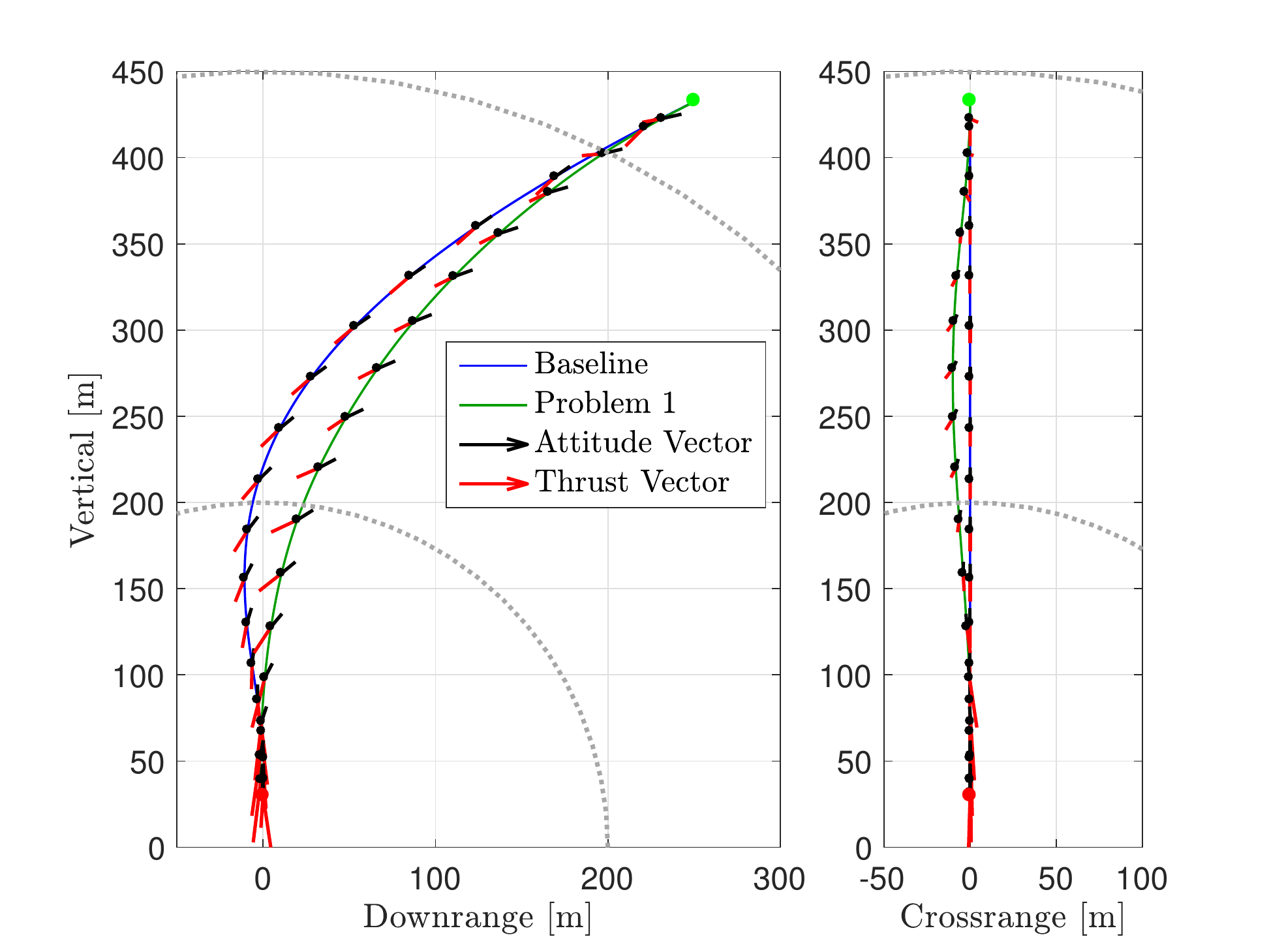}}\hfil
    \subfloat[]{\includegraphics[width=0.48\textwidth]{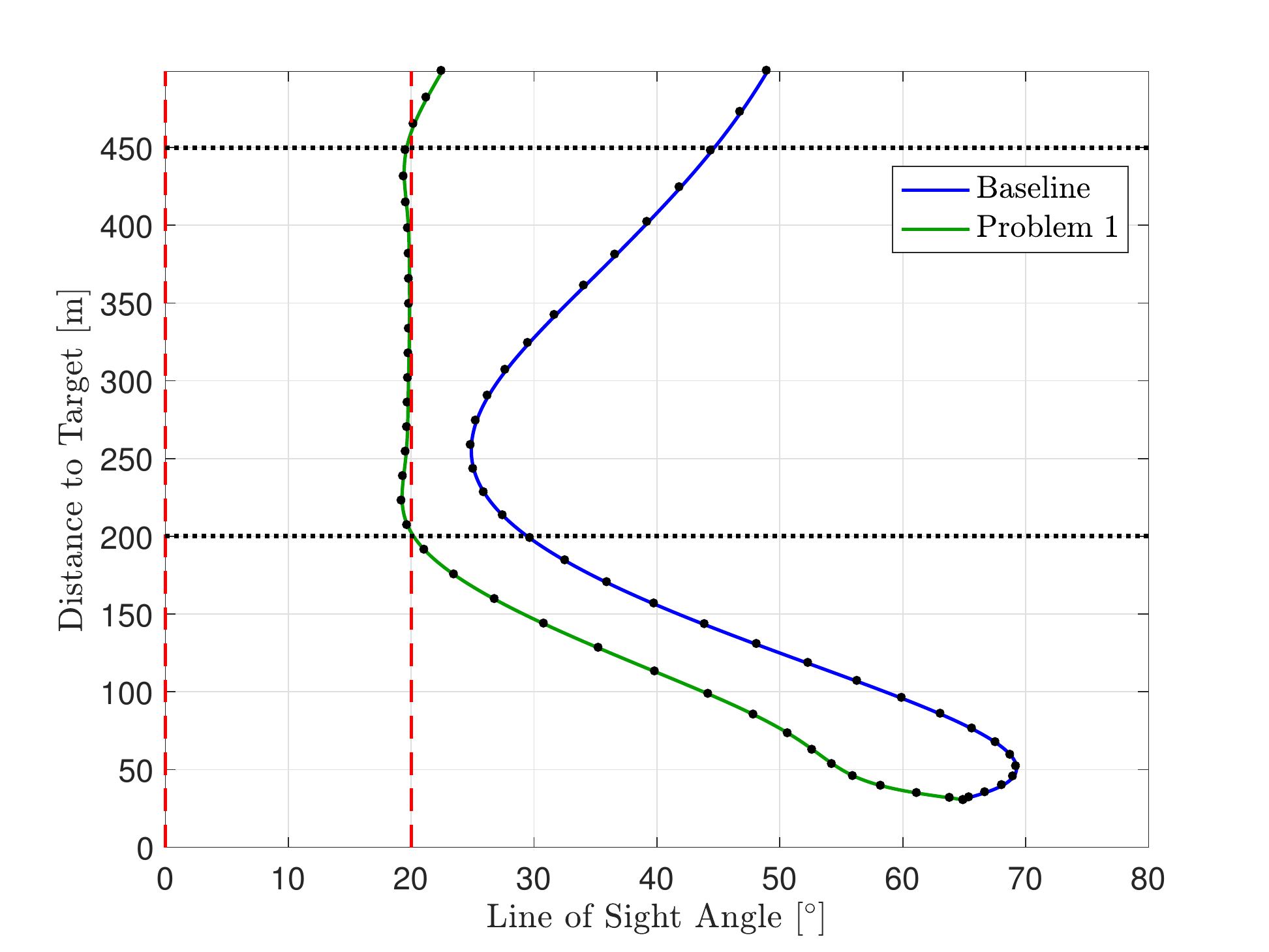}}
    \caption{Trajectories for the baseline problem and Problem~\ref{prob:prob1} with associated attitude and thrust directions. Only half of the discrete points are shown in the right-hand axes for clarity. Also shown in the right-hand axes are the circular boundaries in which the slant-range-triggered constraint is enforced. The time history of the \LOS angle versus slant range is shown in left-hand axes.}
    \label{fig:stc_3d_trj}
\end{figure}

\begin{figure}[tb]
    \centering
    \subfloat[]{\includegraphics[width=0.35\textwidth]{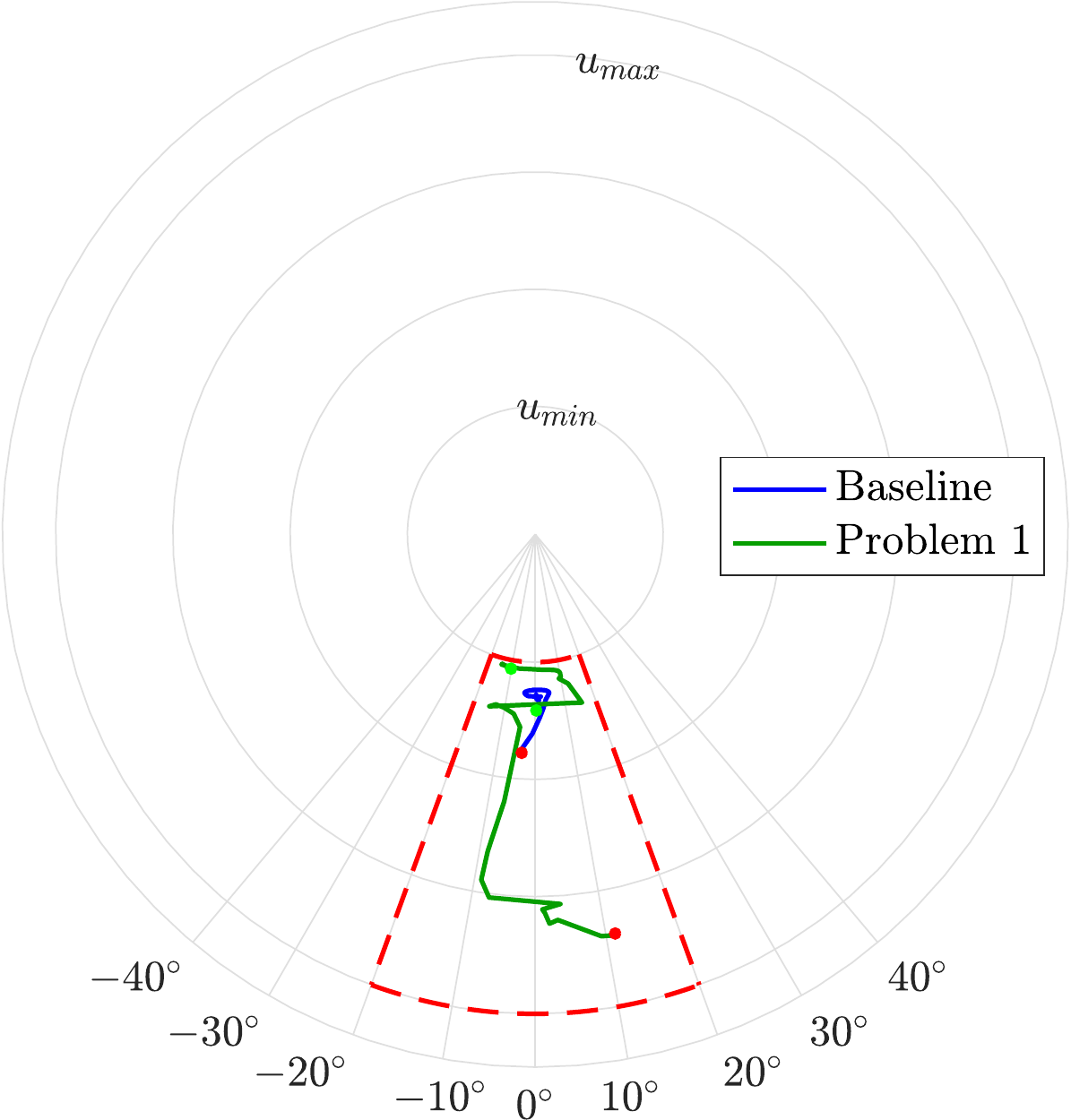}}\hfil
    \subfloat[]{\includegraphics[width=0.5\textwidth]{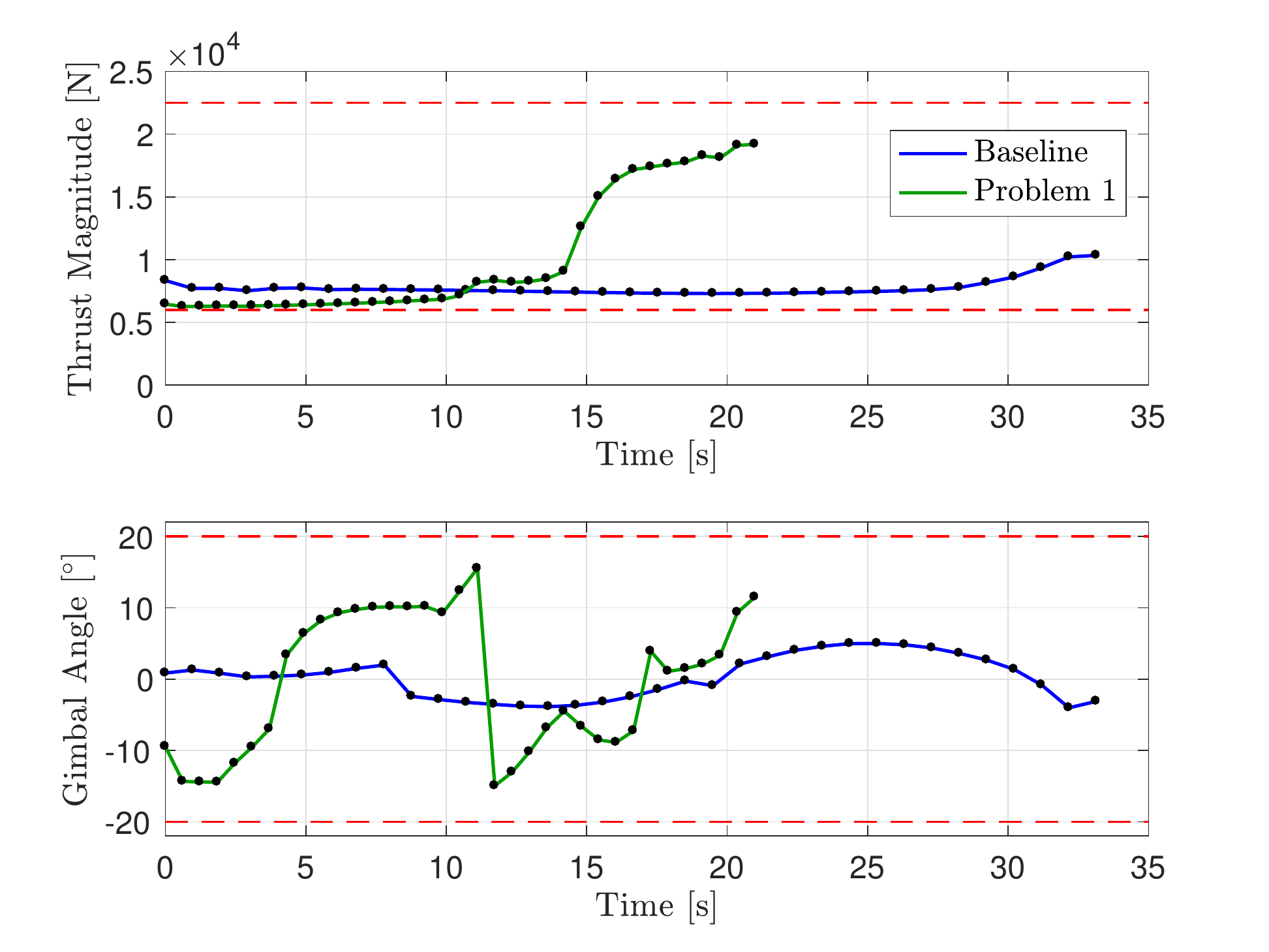}}
    \caption{Polar plot of thrust magnitude versus gimbal angle and their associated time histories for the solutions to the baseline problem and Problem~\ref{prob:prob1}.}
    \label{fig:stc_thrust}
\end{figure}

In this case study we solve both the baseline problem and Problem~\ref{prob:prob1} for comparison. Figure~\ref{fig:stc_3d_trj} shows the resulting trajectories using a temporal density of $N=35$. The black dots represent the discrete solution obtained from the final solve step (i.e., the last solution of Problem~\ref{prob:prob2}), while the solid curves are the result of integrating the discrete controls through the nonlinear dynamics. Both solid curves pass through each discrete point, indicating that the final propagation step has exactly captured the nonlinear dynamics of the problem. Figure~\ref{fig:stc_3d_trj}(b) shows the \LOS angle versus slant range and shows that the optical sensor cannot view the landing site throughout the entire baseline maneuver. In contrast, the solution to Problem~\ref{prob:prob1} satisfies the \LOS constraint between $450\unit{m}$ and $200\unit{m}$ as desired. A deviation from the baseline trajectory is required to ensure feasibility with respect to this constraint, and the differences between the inertial trajectories can be seen clearly in Fig.~\ref{fig:stc_3d_trj}(a). Note that the state-triggered constraint induces some crossrange motion that is otherwise not observed during the planar baseline maneuver. The corresponding thrust curves can be seen from Fig.~\ref{fig:stc_thrust}, which provides both a polar plot of thrust magnitude versus gimbal angle and the individual time histories. Problem~1 uses near-minimum thrust at the beginning of the maneuver, while using the available gimbal to maintain a feasible attitude with respect to the pointing constraint. As soon as the state-triggered constraint is inactive, there is a small divert to achieve a successful landing, observed by a change in the gimbal angle and increased throttle to slow the vehicle for a soft landing.

\begin{figure}[tb]
    \centering
    \subfloat[Baseline Problem]{\includegraphics[width=0.48\textwidth]{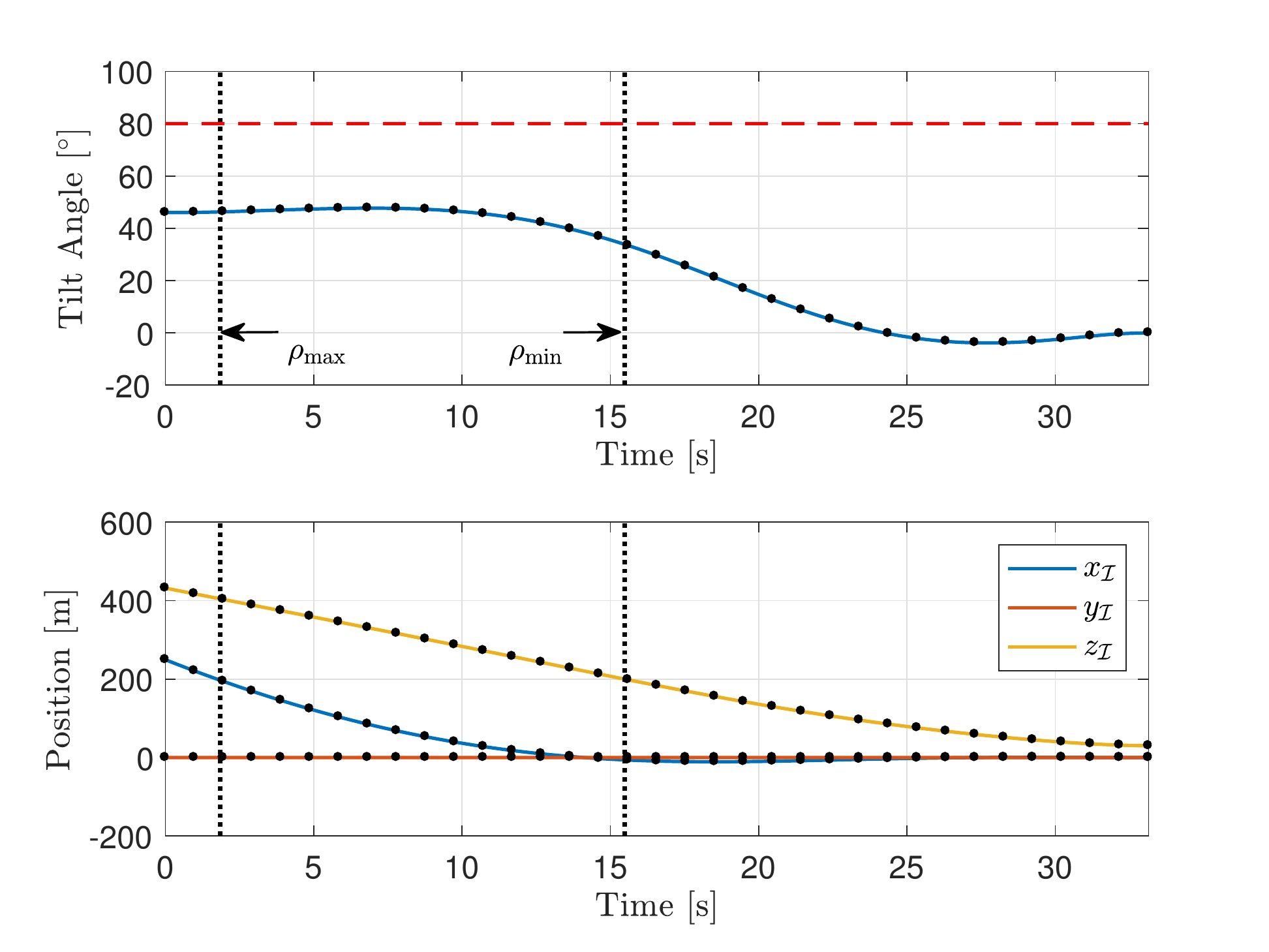}}\hfil
    \subfloat[Problem 1]{\includegraphics[width=0.48\textwidth]{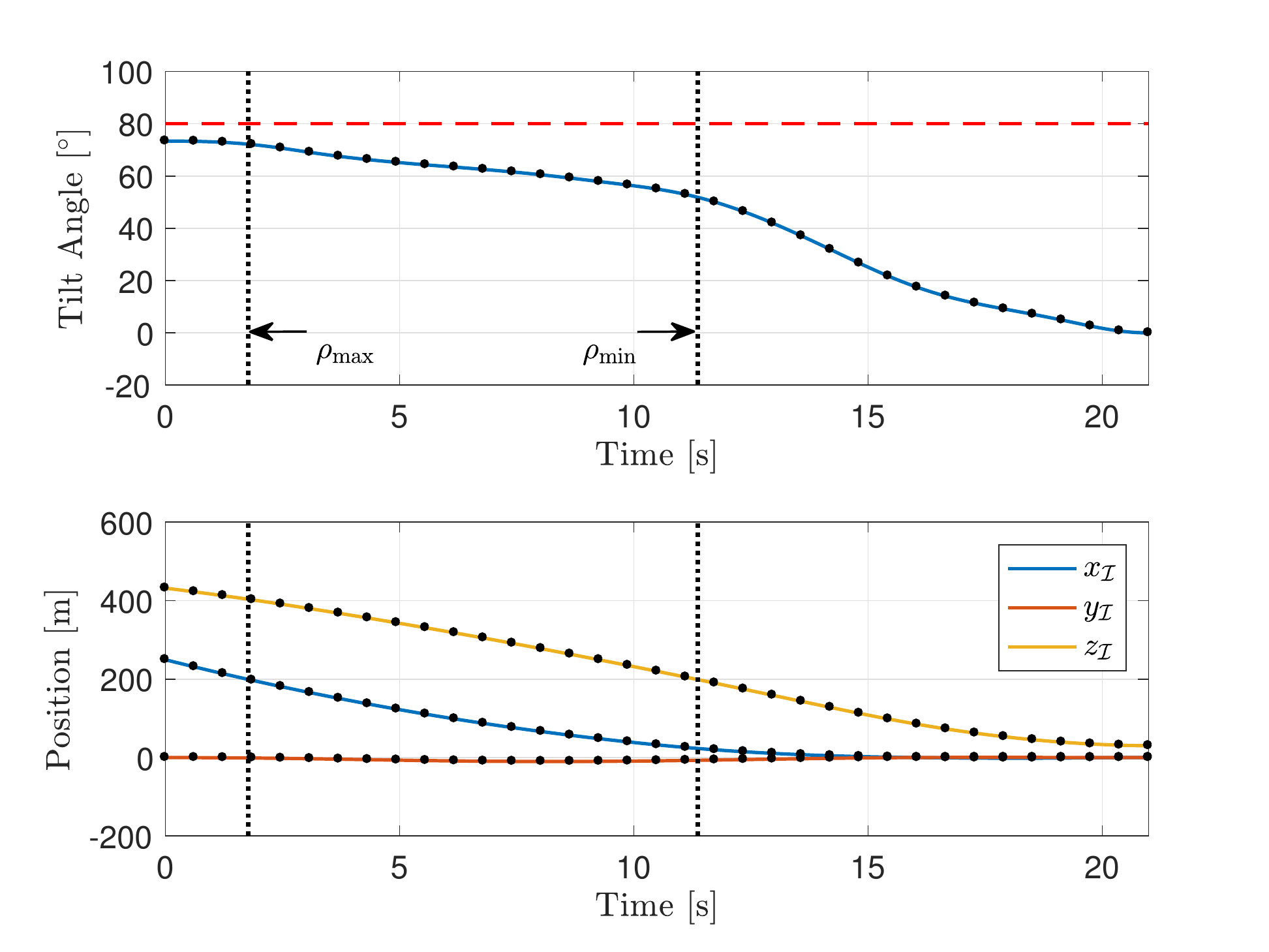}}
    \caption{Vehicle tilt angle and inertial position time histories for solutions to the baseline problem and Problem~\ref{prob:prob1}.}
    \label{fig:stc_states}
\end{figure}

Figure~\ref{fig:stc_states} provides the time history of the vehicle's tilt angle (see~\eqref{eq:tilt_max} and~\eqref{eq:dq_tilt_max}) and inertial position for each maneuver. One can see that when the slant-range-triggered constraint is enforced in Problem~\ref{prob:prob1}, the time spent further than $200\,\unit{m}$ from the landing site is reduced by several seconds. Interestingly, the burn time for the solution to Problem~\ref{prob:prob1} is much lower that the baseline case. 

\subsection{Monte Carlo Case Study}\label{sec2:monte_carlo}

Our objective in this section is to study the ability of the algorithm to successfully compute trajectories given a wide range of initial conditions. At the same time, we aim to assess the real-time capabilities by analyzing the total algorithm runtime. The baseline problem is used for this case study. The optimal trajectory -- and hence algorithm performance -- depends solely on the initial conditions for a given solver, discretization method, initial trajectory guess method, set of algorithm weights, terminal conditions and temporal density. We therefore investigate perturbations from the nominal initial conditions in Table~\ref{tab:params_ex} of sufficient size to test the algorithms capabilities. We use the term \textit{trial} to refer to the converged solution for a single initial condition (i.e., each run of the Monte Carlo study is one trial).

Since the attitude is optimized as part of the solution, we do not consider perturbations to $\q(\ti)$. The initial angular velocity is also kept constant for each trial, since we assume that it may always be controlled to near zero during the preceding descent stage. The initial mass, position, and velocity are different for each trial. The mass is assumed to vary according to a uniform distribution between $90\%$ and $110\%$ of its nominal value
\begin{equation}
    m(\ti) = \mi ( 1+\delta m) \quad \text{where} \quad  \delta m \sim \mathcal{U}(-0.1,0.1).
    \label{eq:mass_var}
\end{equation}
Next, we use the method described in~\cite{Malyuta2019} to compute the initial inertial position and velocity. The inertial velocity is chosen first according to 
\begin{equation}
    \vI(\ti) = \vi{\inertial} + \delta \bm{v},
    \label{eq:vel_var}
\end{equation}
where the dispersions are assumed to come from a zero-mean normal distribution with independent standard deviations of $7\,\unit{m/s}$ in the horizontal directions and $4\,\unit{m/s}$ in the vertical direction. The inertial position $\rI(\ti)$ is then generated via hit-and-run sampling~\cite{Vempala2008} of the constrained controllability polytope~\cite{Dueri2014a,Eren2015} generated for the simplified 3-DoF problem~\cite{Acikmese2007a}. This amounts to sampling a random position from the feasible set of the 3-DoF problem~\cite{Acikmese2007a} that corresponds to the velocity from~\eqref{eq:vel_var}, initial mass from~\eqref{eq:mass_var}, and control constraints in~\eqref{eq:thrust_bound}.

\begin{figure}[tbh]
    \centering
    \subfloat[Open-loop position error.]{\includegraphics[width=0.33\textwidth]{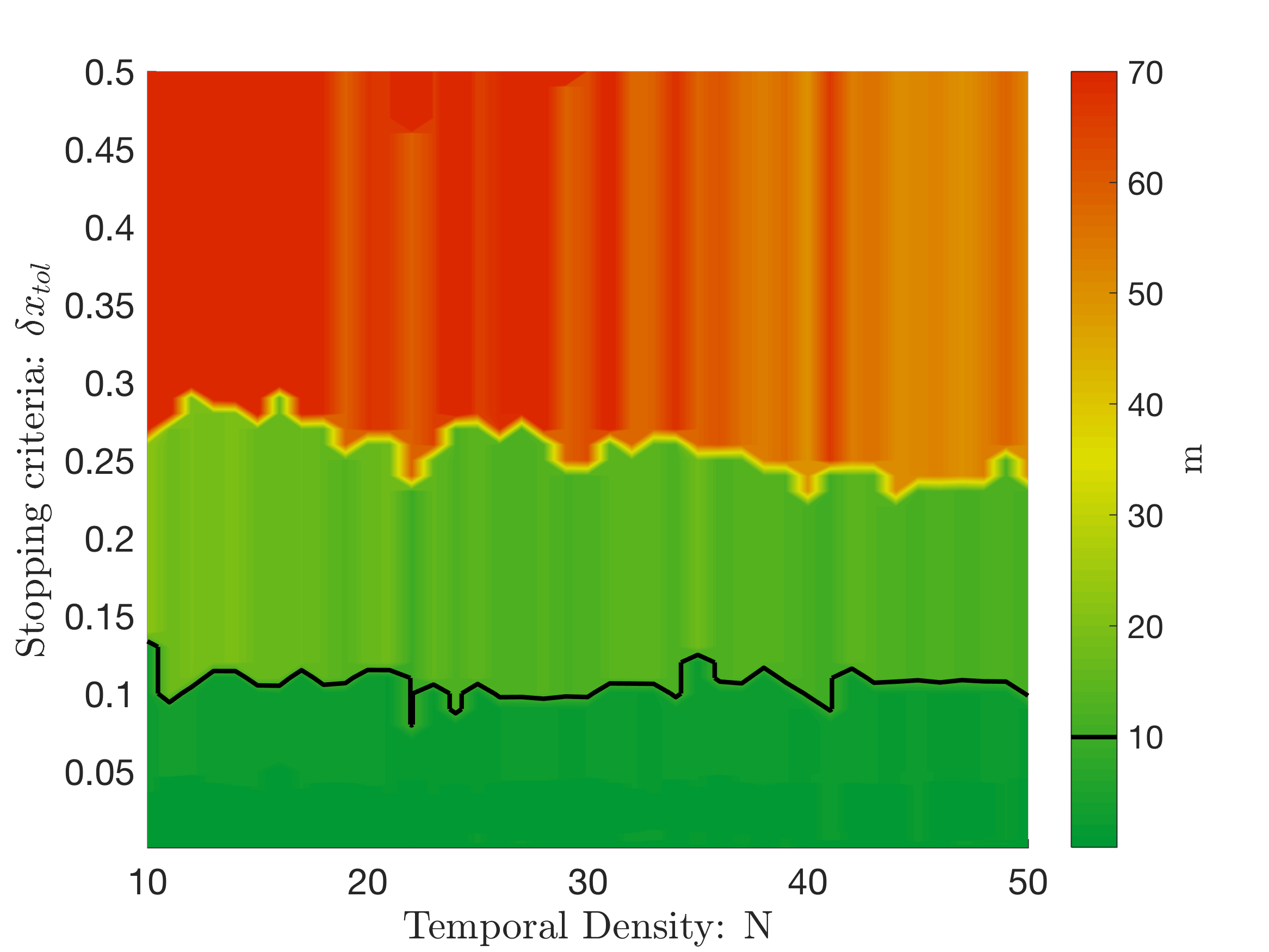}} \hfil
    \subfloat[Open-loop velocity error.]{\includegraphics[width=0.33\textwidth]{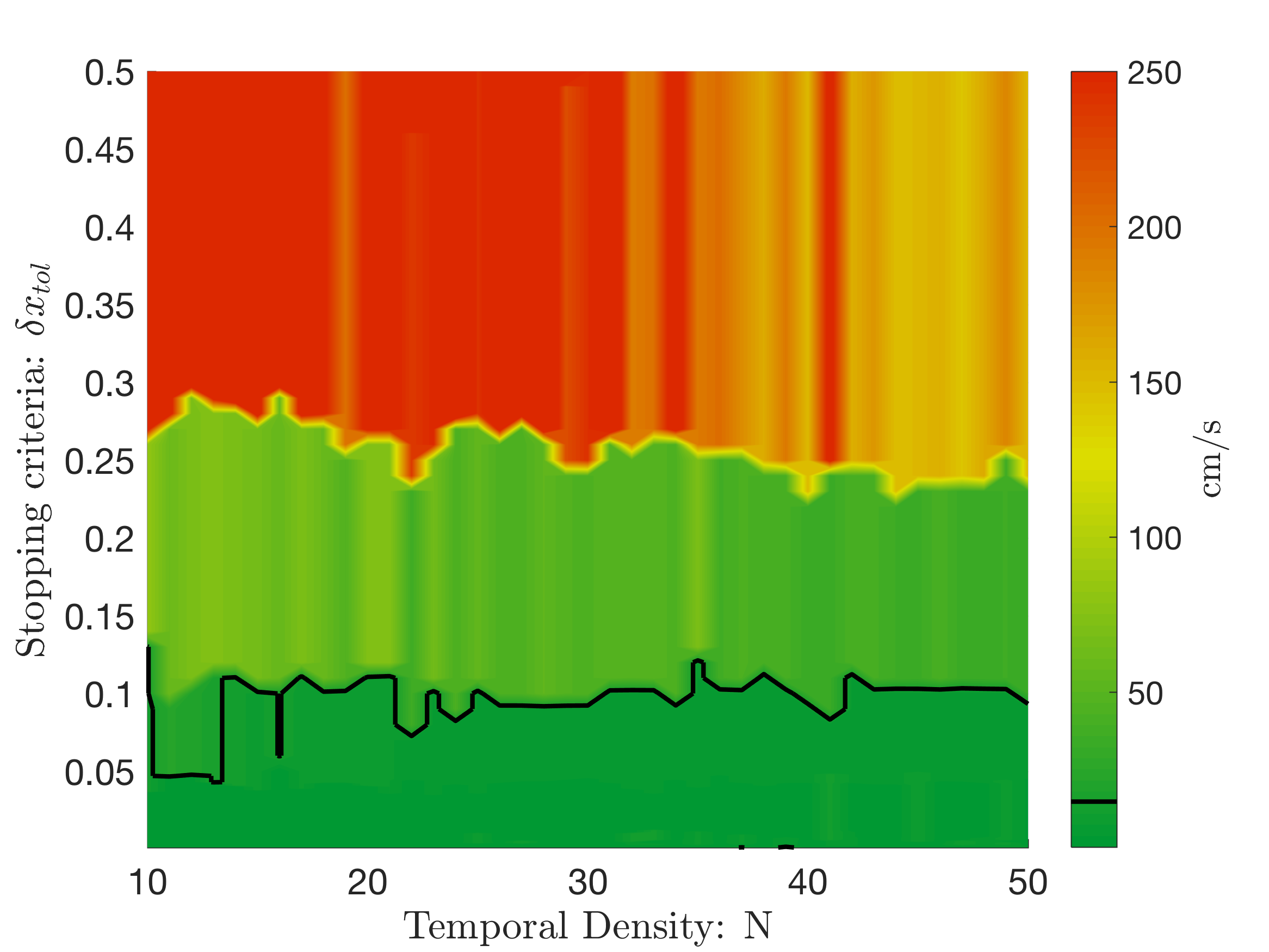}}\hfil
    \subfloat[Total solution time.]{\includegraphics[width=0.33\textwidth]{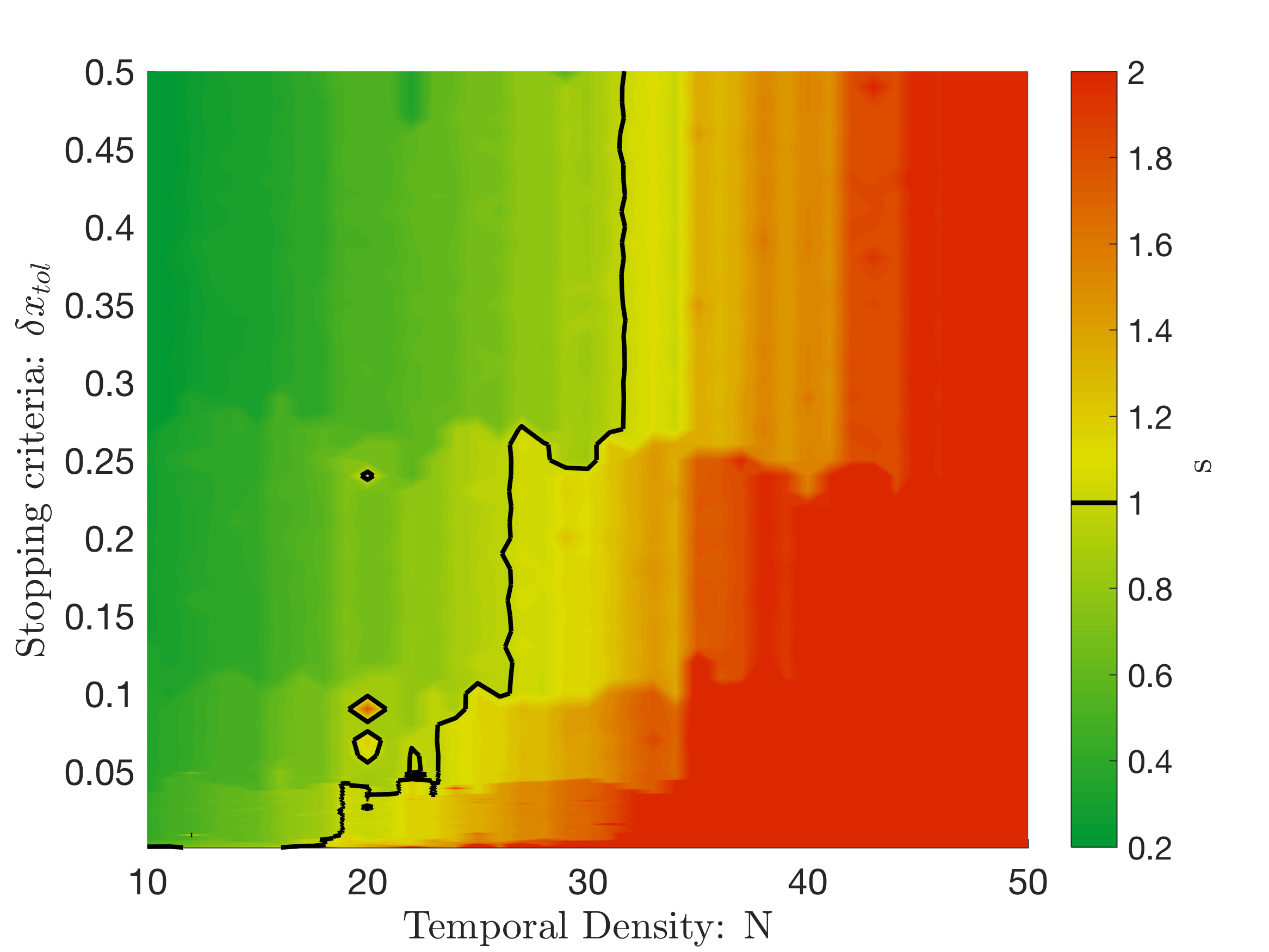}}
    \caption{Open-loop position and velocity errors as a function of temporal density $N$ and stopping criteria $\dxtol$.}
    \label{fig:3d_plots}
\end{figure}

We first investigate the relationship between convergence tolerance, $\dxtol$, the temporal density, $N$, and the resulting open-loop error in the position and velocity states. For each trial, once a converged solution is obtained, the thrust commands are integrated through the original nonlinear dynamics~\eqref{eq:mass_dyn},~\eqref{eq:dq_kinematics} and~\eqref{eq:dq_dynamics}, and the error is defined as the difference between the (integrated) final state and the desired boundary condition. In the limit as $N \rightarrow \infty$ and $\dxtol \rightarrow 0$, the open-loop error is expected to tend to zero at the expense of computation time. Fig.~\ref{fig:3d_plots} shows the open-loop position error, open-loop velocity error and the total CPU runtime as functions of $N$ and $\dxtol$. These plots were generated for each $N$ between $10$ and $50$, while $\dxtol$ was varied in step-sizes of $0.001$ from $0.001$ to $0.05$, and then in step-sizes of $0.01$ between $0.06$ and $0.5$. 

Figure~\ref{fig:3d_plots} can be used to choose key algorithm parameters. Given a set of open-loop error requirements, such as
\begin{equation}
    \| \rI(\tf) - \rf{\inertial} \| \leq 10\,\unit{m}\,, \qquad \| \vI(\tf) - \vf{\inertial} \| \leq 15\,\unit{cm/s}\,,
    \label{eq:success}
\end{equation}
one may select a corresponding stopping criteria $\dxtol$ from Fig.~\ref{fig:3d_plots}(a) and Fig.~\ref{fig:3d_plots}(b). It can be seen that the error bounds in~\eqref{eq:success} can be achieved across a variety of temporal densities so long as $\dxtol\leq0.05$. The \textit{smallest} temporal density that can achieve a given computation time requirement may then be selected. For a computation time of less than $1$ second, Fig.~\ref{fig:3d_plots}(c) shows that it is sufficient to choose $N$ roughly between $10$ and $18$. While precise choices may require more in-depth analysis -- such as the possibility for constraint violation between discrete time nodes -- we proceed with $\dxtol=0.01$ and $N=10$ in what follows.

%

\subsubsection{Dispersed Performance}\label{sec3:dispersed_performance}

\begin{figure}[tbh]
    \centering
    \subfloat[Initial condition dispersions with respect to nominal. ]{\includegraphics[width=0.48\textwidth]{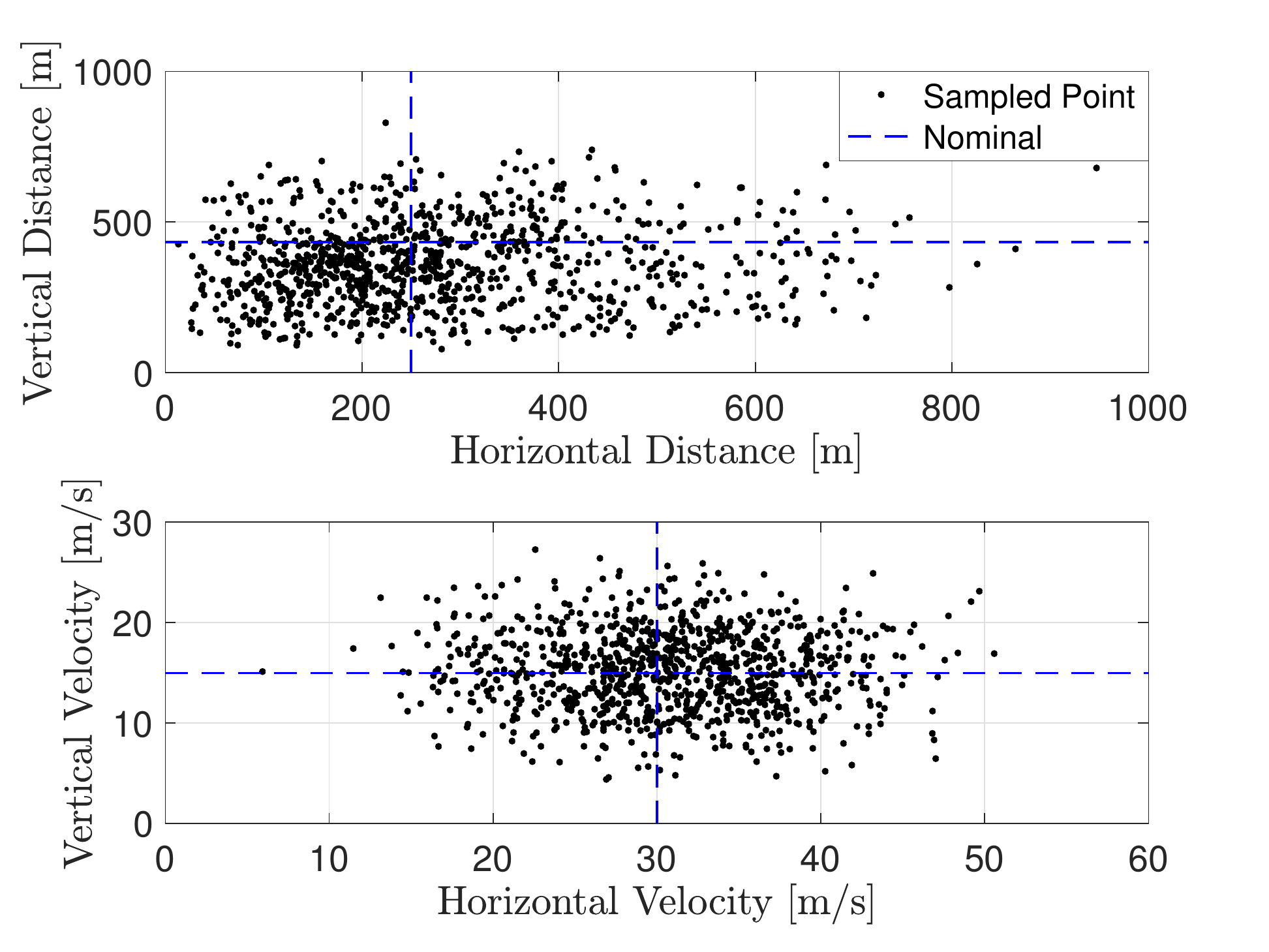}}\hfil
    \subfloat[Final position and velocity errors relative to open-loop objectives.]{\includegraphics[width=0.48\textwidth]{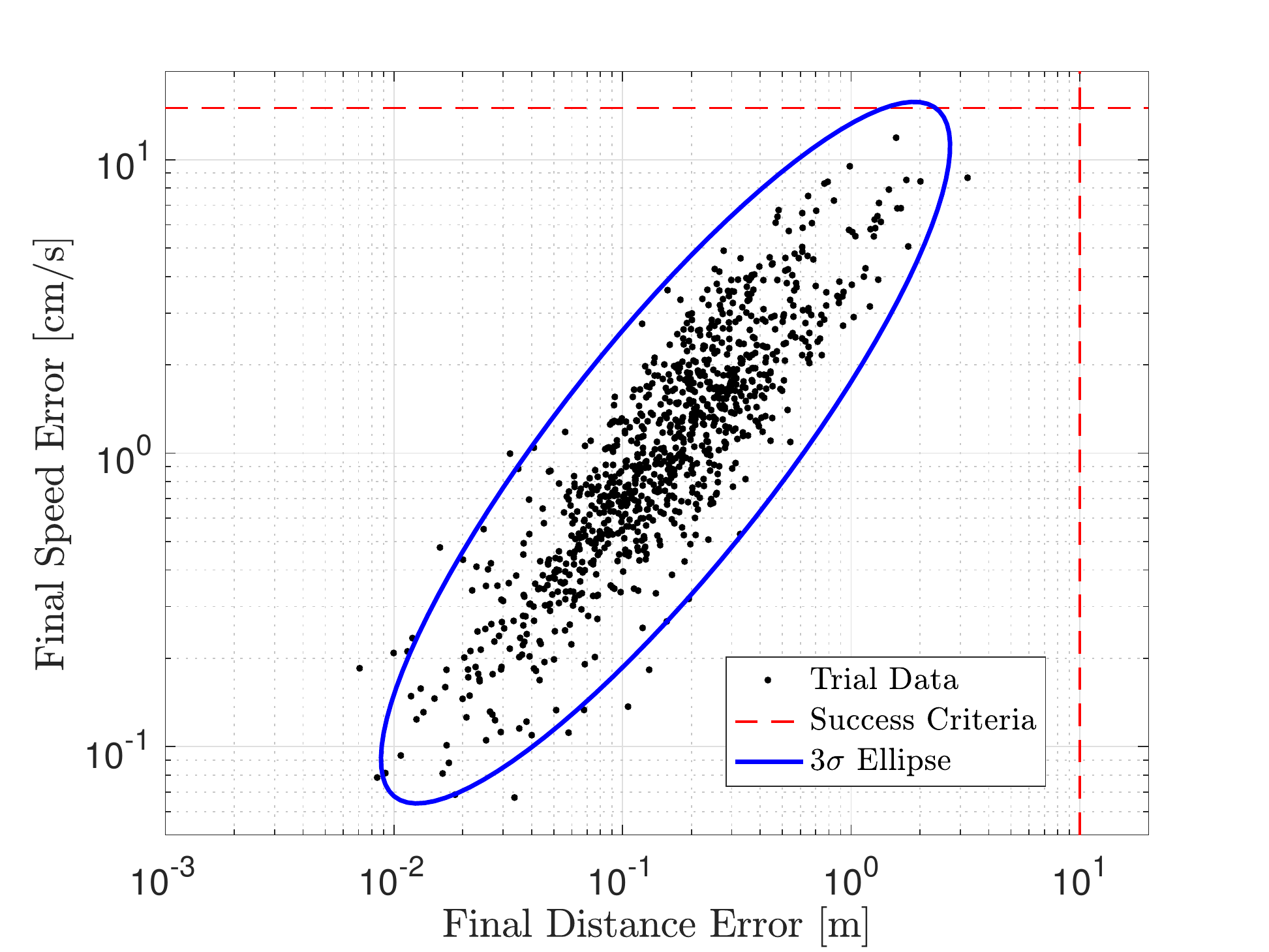}}
    \caption{Position and velocity values at the beginning and end of each successful trial.}
    \label{fig:results_3sig_box}
\end{figure}

We ran 1000 trials using the initial mass, position and velocity dispersions noted above. If a trial met the open-loop requirements in~\eqref{eq:success} it is labeled as ``success'', otherwise it is labeled as ``failed''. Figure~\ref{fig:results_3sig_box}(a) shows the dispersed initial positions and velocities computed using our sampling method discussed above. With respect to the objectives stated in~\eqref{eq:success}, Fig.~\ref{fig:results_3sig_box}(b) shows the final open-loop error in position and velocity for each successful trial. Data analysis reveals that the open-loop position and velocity errors follow a lognormal distribution. We therefore compute statistics on the (natural) logarithm of these errors and shall use these to discuss how well the results achieve our objectives. Figure~\ref{fig:results_3sig_box}(b) shows the three standard deviation ($3\sigma$) confidence ellipse compared to the success criteria~\eqref{eq:success}. One can see that there is a slight overlap with the velocity requirement, indicating that the open-loop objectives are met with slightly less than $99.7\%$ confidence. 

Of the 1000 trials, 971 were successful based on~\eqref{eq:success}. Given that the interpolated initial guess is very simple and generally quite far from being feasible, this is not overly surprising. By taking the 29 failed trajectories and re-solving each of them using the 3-DoF initial guess method~\cite{SzmukReynolds2018}, success was achieved in 28 of these previously failed cases. A perfect success rate of the 3-DoF initialized algorithm was reported in~\cite{Malyuta2019}, albeit for an implementation with a more robust solver, SDPT3, that is not suited to real-time applications. Our results have shown here that our real-time implementation can achieve a 99.9\% success using the same methods, and $97.1\%$ for a computationally simpler method.

\subsubsection{Computational Performance}\label{sec3:computation_performance}

\begin{figure}[tbh]
    \centering
    \subfloat[Total solution time for all dispersed trajectories. Each bin width is $0.1$ seconds.]{\includegraphics[width=0.48\textwidth]{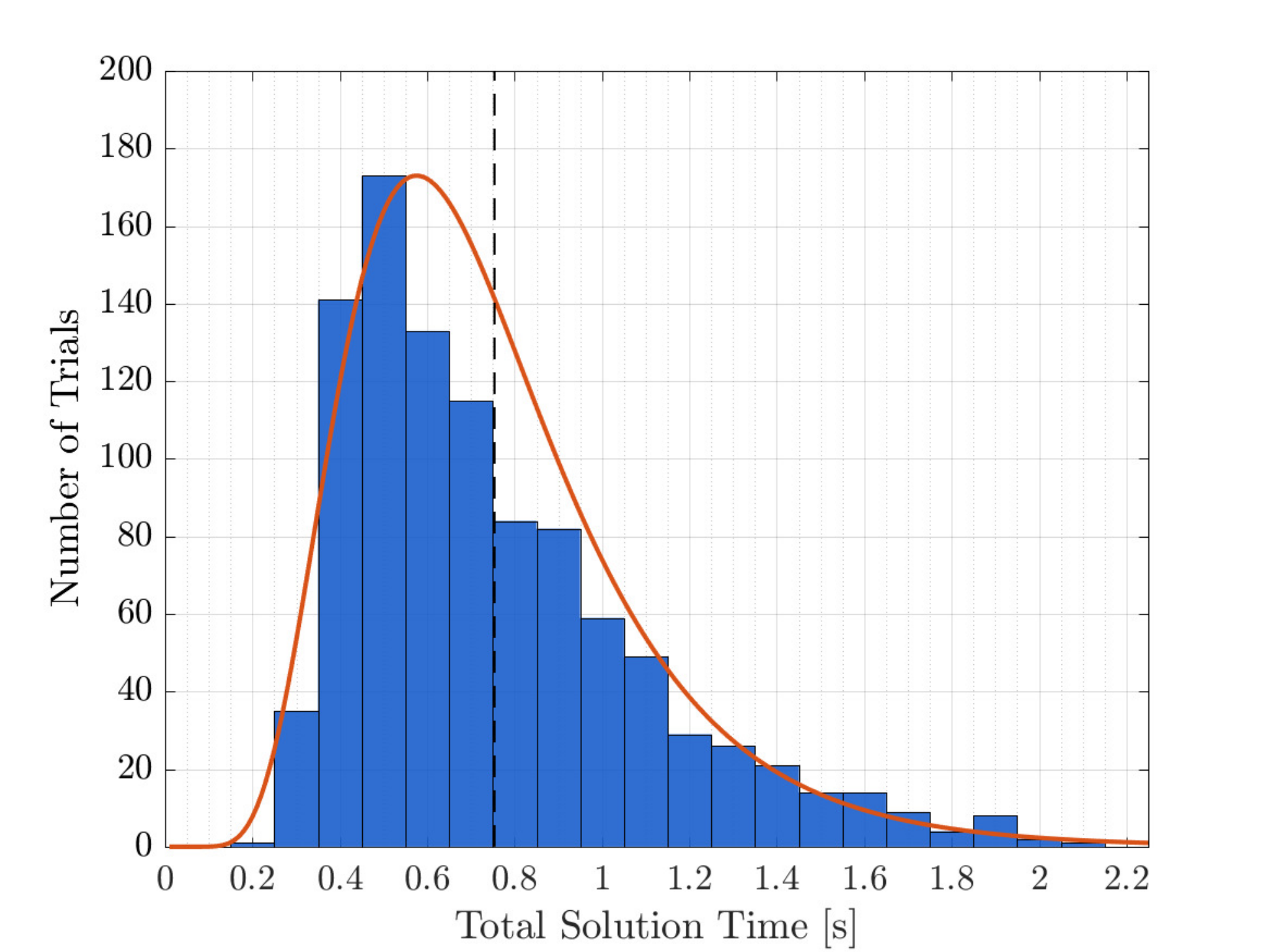}}\hfil
    \subfloat[Q-Q plot of the log of total solution time showing a platykurtic normal distribution.]{\includegraphics[width=0.48\textwidth]{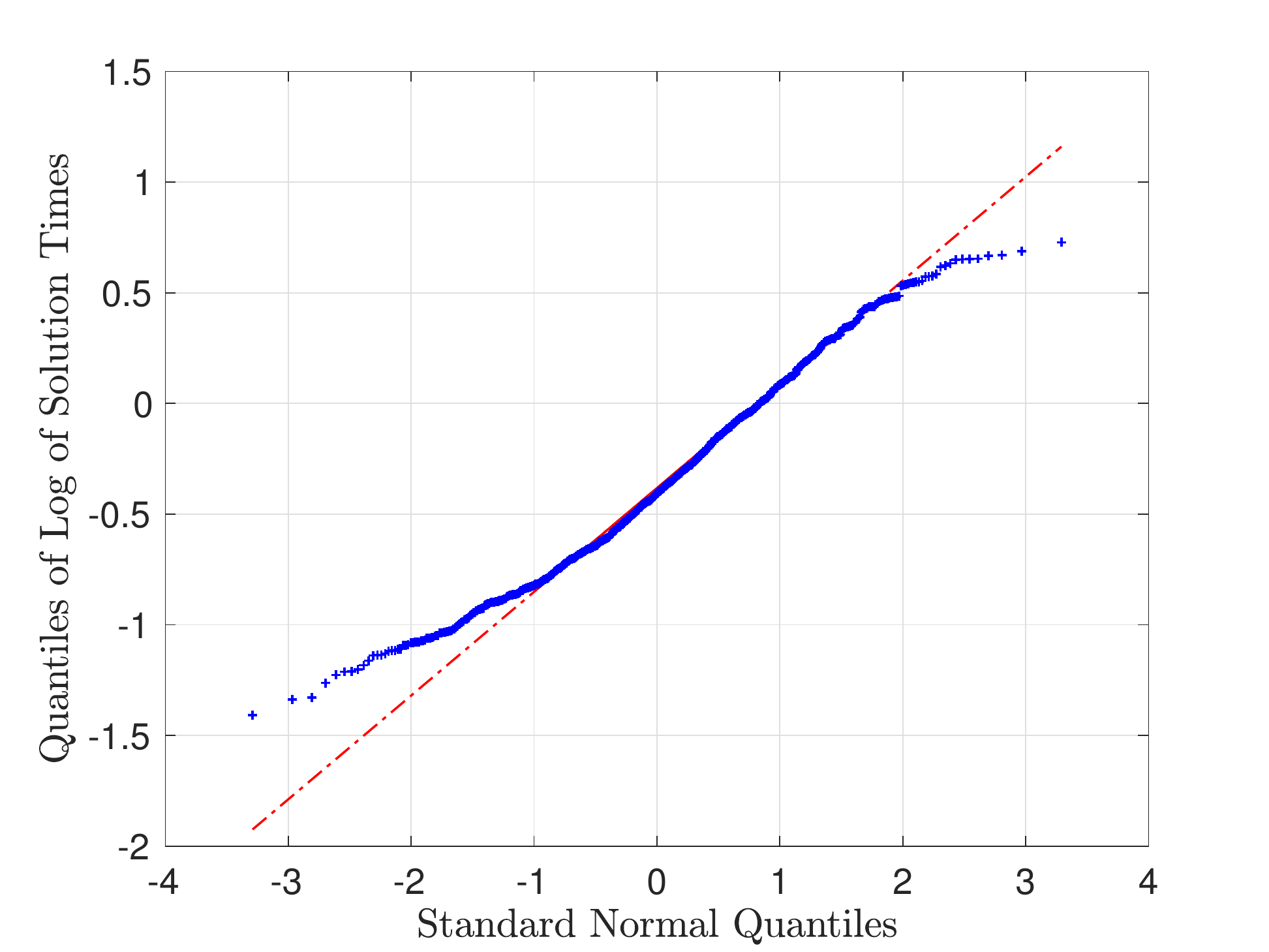}}
    \caption{Computational performance results for all successful trials.}
    \label{fig:results_iters_ttime}
\end{figure}

The results in these tests were all generated on a 2015 iMac with a 3.2 GHz Intel Core i5 processor with 8 GB of RAM. Each trial was executed using a single-core C++ implementation of the algorithm with a MEX-interface to MATLAB\textsuperscript{{\textregistered}}. Figure~\ref{fig:results_iters_ttime}(a) shows the distribution of the total solution time across all trials. The total time is computed as the sum of the propagation and solve step CPU times over all iterations for a single trial. The quantile-quantile (Q-Q) plot in Fig.~\ref{fig:results_iters_ttime}(b) shows that the total solution time is lognormally distributed.  Using this distribution, the mean solution time was computed to be $0.75$ seconds (shown in Fig.~\ref{fig:results_iters_ttime}(b)) with $3\sigma$ confidence that the solve time is less than $2.21$ seconds from a one-tailed analysis.  Moreover, the S-shape of the Q-Q plot indicates \textit{platykurtic} behaviour (i.e., compared to a normal distribution, this distribution produces fewer and less extreme outliers). This behaviour is desirable for flight algorithms as it reduces the probability that excessively long solve times will be observed.

To highlight the contributions of the solve and propagation steps to the total solution time, statistics were computed using all iterations across all trials. The propagation step was executed, on average, in $8.3\pm0.9$ milliseconds, while the solve step took $64.8\pm15.3$ milliseconds. This indicates that the solve step averaged $89\%$ of the total iteration time in our trials. We note that parallel computing techniques can be readily applied to the propagation step, a process that could therefore result in an approximate $10\%$ drop in total solution time. 

\section{Conclusions}\label{sec:conclusion}

This paper has presented a formulation of the 6-DoF powered descent guidance problem using dual quaternions. Both a baseline set of state and control constraints and a novel slant-range-triggered line of sight constraint were introduced. The resulting non-convex free-final time optimal control problem was solved by using an iterative procedure that transcribed the problem into a sequence of convex relaxations. A scaling procedure was introduced to maintain proper numerical conditioning for a wide range of conditions, and a new heuristic method was presented that has been found to guide the iterative process quickly towards convergence. Two numerical case studies looked at the capabilities of the algorithm for a lunar descent scenario. It was shown how a slant-range-triggered constraint can alter a trajectory in both position and attitude so that a line of sight constraint is satisfied in the desired slant range interval. Second, a Monte Carlo analysis established the connection between open-loop error requirements, the algorithm's stopping criterion, and temporal density. The algorithm was tested over a wide range of initial conditions, and was able to meet the desired error bounds of $10~\unit{m}$ and $15~\unit{cm/s}$ with $3\sigma$ confidence. An analysis of the associated computation time shows a mean solution of $0.75~\unit{s}$ on a desktop computer. 

Future work will have two primary directions. The first will focus on the continued study of the optimality, robustness, and convergence properties of the algorithm presented herein. The second will focus on developing a preliminary real-time C implementation that can be used to more accurately characterize the capabilities of the algorithm on representative spacecraft hardware. This implementation will play a key role in reaching accurate and meaningful conclusions about the real-time capabilities of the proposed methodology.

\section*{Acknowledgements}

This research has been supported by NASA grant NNX17AH02A and government sponsorship is acknowledged. A portion of this research was carried out at the Johnson Space Center and the first author would like to thank the members of the Aeroscience and Flight Mechanics branch for their helpful discussions.

\end{document}